\documentclass{amsart}
\usepackage{amsmath}
\usepackage{amssymb}
\usepackage{amsthm}
\usepackage{amstext}
\usepackage{a4wide} 
\usepackage[utf8]{inputenc}
\usepackage{mathtools}

\usepackage{hyperref}
\usepackage{cleveref}
\usepackage{bbm}
\usepackage{xcolor}
\usepackage{mathrsfs}

\usepackage{verbatim}
\usepackage{graphicx}
\usepackage{subcaption}
\usepackage{algorithm}
\usepackage{algorithmic}
\usepackage[colorinlistoftodos]{todonotes}
\usepackage{xpatch}
\newcommand{\T}{\mathrm{T}}       %
\newcommand{\dt}{h}               %

\def\LIN{\mathrm{Lin}}

\newcommand{\mL}{\mathcal{L}}  %
\newcommand{\mY}{\mathcal{Y}}  %
\newcommand{\mS}{\mathcal{S}}  %
\newcommand{\mA}{\TAOP}
\newcommand{\mKS}{\Psi_s}       %

\newcommand{\bE}{\mathbb{E}}   %
\newcommand{\bP}{\mathbb{P}}   %
\def\R{\mathbb{R}}               %
\def\LOPERH{\mathcal{L}_{\dt}}   %
\def\LOPERHPI{\widetilde{\mathcal{L}}_{\dt}}  %
\def\LOPERTIL{\bar{\mathcal{L}}} %
\def\DM{{\Delta \mathcal{M}}}                %
\def\DN{{\Delta \mathcal{N}}}                %
\def\DMN{{\Delta \mathcal{M}_{n}}}

\def\DMM{{\Delta \mathcal{M}_{m}}}

\def\BIGO{\mathcal{O}}          %
\def\BIGOM{\mathcal{O}_P}
\def\WHF{\widehat{f}_{\dt}}
\def\WTF{\widetilde{f}_{\dt}}

\def\WTG{\bar{g}}       %

\def\DMNFH{{\Delta \mathcal{M}_{n}(\WHF)}}   %
\def\DMNFT{{\Delta \mathcal{M}_{n}(\WTF)}}   %
\def\LGC{\left(g-\frac{\dt}{2}\mL g\right)}

\def\TBOP{{\widetilde{\mathcal{B}}_{\dt}}}
\def\TAOP{{\widetilde{\mathcal{A}}}}
\def\AOP{\mathcal{A}}
\def\TQOP{\widetilde{Q}_{\dt}}
\def\COMMA{{\,,}}
\def\PERIOD{{\,.}}
\def\VPAR{\omega}
\def\NSMPL{N_{\mathrm{smpl}}}
\def\TFIN{T_{\mathrm{fin}}}

\renewcommand{\leq}{\leqslant}
\renewcommand{\geq}{\geqslant}
\newcommand{\dps}{\displaystyle}
\newcommand{\cN}{\mathcal{N}}
\newcommand{\Id}{\mathrm{Id}}

\ifx\IMAJNA\undefined
\newtheorem{theorem}{Theorem}[section]
\newtheorem{corollary}{Corollary}[section]
\newtheorem{proposition}{Proposition}[section]
\newtheorem{remark}{Remark}[section]
\newtheorem{lemma}{Lemma}[section]
\newtheorem{example}{Example}[section]
\fi

\newtheorem{assumption}{Assumption}[section]

\begin{document}
\ifx\IMAJNA\undefined
\title[Martingale product estimators for Langevin dynamics]%
{Martingale product estimators for sensitivity analysis in computational statistical physics}

\author[P. Plech\'a\v{c}]{Petr Plech\'a\v{c}}
\address{Department of Mathematical Sciences, University of Delaware, Newark, DE 19716, USA}
\email{plechac@udel.edu}

\author[G. Stoltz]{Gabriel Stoltz}
\address{CERMICS, Ecole des Ponts, Marne-la-Vallée, France \& MATHERIALS team-project, Inria Paris, France}
\email{gabriel.stoltz@enpc.fr}

\author[T. Wang]{Ting Wang}
\address{Physical Modeling \& Simulation Branch, CISD, DEVCOM Army Research Laboratory,  Aberdeen Proving Ground, MD 21005, USA}
\email{tingw@udel.edu}

\date{\today}



\begin{abstract}
We introduce a new class of estimators for the linear response of steady states of stochastic dynamics.
    We generalize the likelihood ratio approach and formulate the linear response as a product of two martingales,
    hence the name "martingale product estimators".
    We present a systematic derivation of the martingale product estimator, and show how to construct such estimator
    so its bias is consistent with the weak order of the numerical scheme that approximates the underlying stochastic
    differential equation. Motivated by the estimation of transport properties in molecular systems, we present
    a rigorous numerical analysis of the bias and variance for these new estimators in the case of Langevin dynamics.
    We prove that the variance is uniformly bounded in time and derive a specific form of the
estimator for second-order splitting schemes for Langevin dynamics.
For comparison, we also study the bias and variance of a Green-Kubo estimator, motivated, in part, by
its variance growing linearly in time. Presented analysis shows that the new martingale product estimators, having uniformly
bounded variance in time, offer a competitive alternative to the traditional Green-Kubo estimator.
We compare on illustrative numerical tests the new estimators with results obtained by the Green-Kubo method.
\end{abstract}

\maketitle

\else
 \title{Martingale product estimators for sensitivity 
 analysis in computational statistical physics}

\shorttitle{Martingale product estimators for Langevin dynamics} 

 \author{%
 {\sc 
 Petr Plech\'{a}\v{c}\thanks{Email: plechac@udel.edu}}\\[2pt]
 Universpity of Delaware, Newark, DE 19716, USA \\[6pt]
 {\sc and}\\[6pt]
 {\sc
 Gabriel Stoltz\thanks{Email:gabriel.stoltz@enpc.fr}}\\[2pt]
 CERMICS, Ecole des Ponts, Marne-la-Vallée, France \& MATHERIALS team-project, Inria Paris, France\\[6pt]
 {\sc and}\\[6pt]
 {\sc Ting Wang\thanks{Email:tingw@udel.edu}}\\[2pt]
 Physical Modeling \& Simulation Branch, CISD, U.S. Army Research Laboratory,  Aberdeen Proving Ground, MD 21005, USA
 }

\shortauthorlist{P. Plech\'a\v{c}, G. Stoltz, T. Wang}

\maketitle

\begin{abstract}
  {We introduce a new class of estimators for the linear response of steady states of stochastic dynamics.
    We generalize the likelihood ratio approach and formulate the linear response as a product of two martingales,
    hence the name "martingale product estimators".
    We present a systematic derivation of the martingale product estimator, and show how to construct such estimator
    so its bias is consistent with the weak order of the numerical scheme that approximates the underlying stochastic
    differential equation. Motivated by the estimation of transport properties in molecular systems, we present
    a rigorous numerical analysis of the bias and variance for these new estimators in the case of Langevin dynamics.
    We prove that the variance is uniformly bounded in time and derive a specific form of the
estimator for second-order splitting schemes for Langevin dynamics.
For comparison, we also study the bias and variance of a Green-Kubo estimator, motivated, in part, by
its variance growing linearly in time. Presented analysis shows that the new martingale product estimators, having uniformly
bounded variance in time, offer a competitive alternative to the traditional Green-Kubo estimator.
We compare on illustrative numerical tests the new estimators with results obtained by the Green-Kubo method.}
{linear response, Langevin equation, likelihood ratio method, variance reduction, non-equilibrium steady states}
\end{abstract}
\fi

\section{Introduction}

We consider the linear response at infinite time ($t\to\infty$) 
to a perturbation in the drift vector field of an ergodic stochastic dynamics 
defined by a system of stochastic differential equations
\begin{equation}\label{eq:dynamics}
dy_t^\eta = (b(y_t^\eta)+\eta F(y_t^\eta)) dt + \sigma(y_t^\eta) dW_t\,\COMMA
\end{equation}
where~$\eta \in \mathbb{R}$ is a small parameter. Under the assumption that
both the unperturbed ($\eta=0$) and perturbed system are ergodic with respect to probability measures $\pi$ and  $\pi^\eta$, respectively, the central object of the numerical analysis we present is the estimation of the quantity
\begin{equation}\label{eq:linresponse}
\left.\frac{d}{d\eta}\right|_{\eta=0}\pi^{\eta}(f) \equiv \lim_{\eta \to 0} \frac{1}{\eta}\left[\int f(y) \,\pi^\eta(dy) - \int f(y)\,\pi(dy)\right]\,,
\end{equation}
for a suitable class of functions $f$, termed observables. The well posedness of this limit is guaranteed under certain assumptions on the dynamics, see for instance~\cite{HM10,lelievre2016partial}.

Our motivation for developing an efficient numerical estimation of~\eqref{eq:linresponse} is the computation of transport coefficients in molecular systems at thermal equilibrium in the canonical ensemble. Accurately estimating transport properties such as the mobility, the shear stress or 
the thermal conductivity, presents computational challenges in spite of algorithmic developments in computational methods for statistical mechanics of fluids and gasses, see, \emph{e.g.}, \cite{evans-moriss,tuckerman2010statistical} for physical background, \cite[Section~5]{lelievre2016partial} for mathematical background and analysis, and \cite{leimkuhler2016computation} for elements of numerical analysis.
We present in Section~\ref{sec:general-continuous-setting} a strategy for deriving a new type of estimators for the problem
\eqref{eq:linresponse} under the dynamics~\eqref{eq:dynamics}. 
A systematic but formal derivation is followed by rigorous numerical analysis and results  in Sections~\ref{sec:MP_for_Langevin} and~\ref{sec:LR-vs-GK}, focusing specifically on the Langevin dynamics 
\begin{equation}\label{eq:Langevinintro}
\left\{ \begin{aligned}
    dq^\eta_t &= M^{-1} p^\eta_t \, ,\\
    dp^\eta_t &= \Big(-\nabla V(q^\eta_t)+\eta F(q^\eta_t,p^\eta_t)\Big)dt 
              -\gamma M^{-1} p^\eta_t dt + \sqrt{\frac{2\gamma}{\beta}} dW_t\,,
\end{aligned} \right.
\end{equation}
which describes particles with positions~$q$ and momenta~$p$ at a fixed inverse temperature
$\beta$. The particle interaction is defined by the potential $V(q)$ and $M$ denotes
the mass tensor.
The specific form of the magnitude~$\sqrt{2\gamma\beta^{-1}} dW_t$ of the standard Wiener process~$W_t$ ensures that the Boltzmann--Gibbs probability measure $\pi(dq\,dp) = Z^{-1} \mathrm{e}^{-\beta(p^T M^{-1} p/2 + V(q))} \,dq\,dp$ is invariant for the unperturbed system at~$\eta=0$.

The perturbing field $F$ is {\it not} assumed to be the gradient of a potential, hence
the corresponding invariant measure $\pi^\eta$ of the perturbed dynamics is unknown in general, and corresponds 
to what is called a {\it non-equilibrium steady state}. As we shall see in Section~\ref{sec:general-continuous-setting} the proposed estimators do not require 
simulations of the perturbed dynamics. Instead the estimator is based on certain ergodic averages  
on trajectories of the {\it unperturbed dynamics} ($\eta=0$). Let us however emphasize that, although we consider as a reference dynamics a Langevin dynamics where the force derives from a gradient, our results can be straightforwardly extended to the situation when~$-\nabla V(q)$ is replaced by a generic (non-gradient) force field, in which case the reference probability measure~$\pi$ is not known.

Standard methods to estimate transport coefficients are reviewed in \cite[Section~5]{lelievre2016partial}. There are essentially two approaches. The first one is to numerically approximate the linear response, which is however computationally quite expensive since averages with respect to~$\pi^\eta$ have to be computed by time averages over very long times in order for the statistical error to scale as~$\eta^{-2}$, so that the limit in~\eqref{eq:linresponse} is sufficiently well approximated. The second popular option is to rely on Green--Kubo formulas, which rewrite the linear response as some integrated correlation function for the equilibrium dynamics. Properly estimating this time integral is however challenging~\cite{dSOG17,EMB17} because the correlation function is very difficult to reliably estimate for long times since it is a small quantity plagued by a large relative statistical error. For these two standard approaches, the numerical analysis of the bias arising from the timestep discretization of the underlying SDE can be read in~\cite{leimkuhler2016computation,lelievre2016partial}. 

An alternative to the GK estimator relies on the so-called likelihood ratio method, which is based on 
the Girsanov theorem. More precisely, the estimator which is suggested in this approach is constructed by linearizing the Girsanov weight considered on a finite time interval, and passing to the infinite time limit. It can be shown that this indeed provides a consistent
estimator of the linear response~\eqref{eq:linresponse}. 
This approach was employed in~\cite{glynn2019likelihood} for discrete time Markov chains and
in~\cite{wang2019steady} for continuous time Markov chains, resulting in 
likelihood ratio type estimators with a variance uniformly bounded in time. 
The derived low-variance modification of the likelihood ratio (LR) method was further 
extended to the linear response of the dynamics described by stochastic differential equations in~\cite{plechac2019convergence}. 
However, the approach developed in \cite{plechac2019convergence} did not provide 
a systematic way to derive estimators of order higher than~1 when discretizing the underlying dynamics.

Our aim in this work is to generalize the likelihood ratio approach and derive a new class of estimators 
that we call {\it martingale product} estimators (MP). The MP estimators reformulate
the linear response as an expectation of the product of two martingales.
The analysis of the proposed approach also provides a systematic way to obtain
estimators whose time-step bias is consistent with the weak order of the underlying stochastic dynamics. In particular, we present estimators that lead to biases of order~2 in the time-step for weak second order discretizations of Langevin dynamics~\eqref{eq:Langevinintro} based on splitting schemes~\cite{leimkuhler2013rational}. Furthermore, we prove that the variance of the MP estimators
is uniformly bounded in time, as for CLR estimators, and therefore provide a competitive alternative to 
Green-Kubo estimators for which the variance grows linearly in time. We more precisely compare the new MP estimators and the Green-Kubo method in Section~\ref{sec:LR-vs-GK} from a mathematical viewpoint, and in Section~\ref{sec:numerics} from a numerical perspective on toy examples. Let us however immediately emphasize that it is not our aim in this study to obtain strong conclusions on the actual relative numerical performances of these two methods. 

The manuscript is organized as follows. We start by presenting in Section~\ref{sec:general-continuous-setting} the heuristic strategy for constructing MP estimators for general stochastic dynamics. In order to make this strategy precise and rigorous for Langevin dynamics, we review some properties of this dynamics and its discretization by splitting schemes in Section~\ref{sec:LD-application}. We can then perform the numerical analysis of the bias and variance of MP estimators for Langevin dynamics in Section~\ref{sec:MP_for_Langevin} (with some technical results postponed to Appendices~\ref{app:useful_estimates} and~\ref{app:secondorder}). For comparison, we also study the bias and variance of the GK estimators in Section~\ref{sec:LR-vs-GK}. Finally, several numerical benchmarks are provided to demonstrate the consistent between our theoretical results and the numerical simulations in Section~\ref{sec:numerics}.

\section{Formal construction of martingale product estimators}
\label{sec:general-continuous-setting}

We present in this section a general strategy for constructing the new MP estimators for the linear response of stationary measures of stochastic dynamics. We first recall in Section~\ref{sec:LR_stoch_dyn} the mathematical framework for the linear response of stationary measures of general diffusion processes. In particular, we express the response coefficients in terms of solutions to Poisson equations and their computation in terms of an auxiliary martingale~\cite{plechac2019convergence}. We then present in Section~\ref{subsec:formal-derivation} 
a general form of the numerical estimators we consider when the underlying continuous diffusion process is discretized in time. Finally, in order to reduce the bias induced by the time discretization on the estimation of the linear response, the key point in our approach is to construct discrete auxiliary martingales tailored to the weak order of the numerical scheme at hand, as formally discussed in Section~\ref{sec:construction_discrete_martingale} for schemes of both first and second weak order.

\subsection{Linear response of stochastic dynamics}
\label{sec:LR_stoch_dyn}

\subsubsection{Reference dynamics.} Given a filtered probability space $(\Omega, \mathcal{F}, \mathcal{F}_t, \bP)$, we denote by $W_t = (W_{1,t},\dots, W_{k,t})^{\T}$ a $k$-dimensional $\{\mathcal{F}_t\}_{t>0}$-adapted standard Wiener process. We study certain perturbations of a $d$-dimensional time homogeneous stochastic differential equation (SDE) on the state space~$\mY$ (for example, $\mY = \mathbb{R}^d$ or~$\mathbb{T}^d$ with $\mathbb{T}=\mathbb{R}\backslash \mathbb{Z}$ the one-dimensional torus):
\begin{equation}\label{eqn:SDE}
dy_t = b(y_t)\, dt + \sigma(y_t) \, dW_t\COMMA
\end{equation}
where $b: \mY \to \mathbb{R}^d$ is a smooth vector field and $\sigma: \mY \to \mathbb{R}^{d \times k}$ is a matrix valued function.
It is implicitly assumed that $b$ and $\sigma$ satisfy all required conditions so that the above SDE is well posed~\cite{oksendal2013stochastic}. We consider the case of degenerate noise, \emph{i.e.} $k < d$, but under the assumption that the matrix field $\sigma(y) = \left(\sigma_{ij}(y)\right)_{1 \leq i \leq d, 1 \leq j \leq k}$ has full rank~$k$ for all~$y \in \mY$.
Equation~\eqref{eqn:SDE} can be written component-wise as
\[
dy_{i,t} = b_i(y_t)\, dt + \sum_{j=1}^k \sigma_{ij}(y_t) \, dW_{j,t}\COMMA \qquad  i=1,\dots,d \PERIOD
\]
We assume that the dynamics~\eqref{eqn:SDE} is ergodic and that its unique invariant measure~$\pi(dy)$ has a smooth positive density~$\rho(y)$ with respect to the Lebesgue measure (see~\cite{kliemann1987recurrence, bellet2006ergodic} for sufficient conditions). 
To write the expression of the generator associated with~\eqref{eqn:SDE}, we introduce the $d\times d$ diffusion matrix~$S = \sigma \sigma^{\T}$, and denote by~$\nabla^2 f$ the Hessian matrix with entries $\partial^2_{y_i,y_j} f$ for $1 \leq i,j \leq d$. 
The symbol ``$:$'' represents the Frobenius inner product of $d\times d$ matrices, \emph{i.e.}, $S:\nabla^2 f = \sum_{1 \leq i,j \leq d} S_{ij} \partial^2_{y_i,y_j} f$. We emphasize that the diffusion matrix $S(y)$ can be degenerate. With this notation the generator of~\eqref{eqn:SDE} acts on smooth test functions~$f$ as  
\begin{equation}
  \mathcal{L}f = b^{\T} \nabla f + \frac{1}{2} S : \nabla^2 f \PERIOD
\end{equation}
{Here and in the sequel, we adopt the following 
notational convention: all vectors are understood as column vectors,  
the gradient of a scalar function $f$ is represented by a column vector $\nabla f$ 
and thus 
$b^T\nabla f = (\nabla f)^\T b = \sum_i b_i \partial_{y_i} f $. 
The notation $\nabla F$ for the gradient of a vector field $F$ represents
a matrix with entries
$(\nabla F)_{ij} = \partial_{y_i} F_j$, i.e., the $j$th column of~$\nabla F$ is the gradient of the individual component~$F_j$ (understood as a column vector).}

\subsubsection{Perturbations of the reference dynamics.}
We next consider perturbations of the reference dynamics~\eqref{eqn:SDE} obtained by adding to the drift field~$b$ a vector field $F: \mY \to \mathbb{R}^d$, with a forcing magnitude~$\eta \in \mathbb{R}$. The perturbed dynamics then reads
\begin{equation}
  \label{eqn:perturbed-SDE}
  dy_t^{\eta} = (b(y_t^{\eta}) + \eta F(y_t^{\eta}))\, dt + \sigma(y_t^{\eta})\, dW_t\PERIOD
\end{equation}
Its generator can be written as
\[
\mL^{\eta} = \mL + \eta \LOPERTIL\COMMA \qquad \LOPERTIL = F^{\T} \nabla\PERIOD
\]
As for the reference dynamics~\eqref{eqn:SDE}, we assume that the perturbed dynamics~\eqref{eqn:perturbed-SDE} is ergodic and that its unique invariant measure~$\pi^\eta(dy)$ has a positive smooth density~$\rho^\eta(y)$ with respect to the Lebesgue measure. The linear response of a given observable of interest~$f$ is then defined as (provided the limit exists)
\begin{equation}
  \label{eqn:lin-response}
  \LIN (f) = \lim_{\eta \to 0}\frac{1}{\eta}(\pi^{\eta}(f) - \pi(f)) =  \lim_{\eta \to 0}\frac{1}{\eta}\int_{\mY} f(y) (\rho^{\eta}(y) - \rho(y)) \, dy\COMMA
\end{equation}
where for the ease of notation $\pi = \pi^0$ and $\rho=\rho^0$ denote the reference measure and reference density, respectively. We emphasize that the linear response depends both on the forcing~$F$ and the reference dynamics.

We present one possible framework for justifying the well-posedness of the limit~\eqref{eqn:lin-response}, relying on linear response theory (see for instance~\cite{HM10}). The reference functional space is the Hilbert space 
\[
L^2(\pi) = \left\{f : \mY \to \mathbb{R}~\text{measurable}~\left|~ \int_{\mY} |f|^2 \, d\pi < \infty\right.\right\}\PERIOD 
\]
Introducing the projection operator 
\begin{equation}
  \label{eqn:projection-operator}
  \Pi f = f - \pi(f),
\end{equation}
we further define the projected subspace
\[
L_0^2(\pi) = \Pi L^2(\pi) = \left\{ f \in L^2(\pi) ~\middle|~ \pi( f )= 0 \right\},
\]
which consists of all $L^2(\pi)$ functions with zero mean with respect to the reference invariant measure~$\pi$. We denote by $\mL^*$ the adjoint of~$\mL$ on~$L^2(\pi)$, \emph{i.e.}, for any $C^{\infty}$ compactly supported test functions~$f, g$, 
\[
\int_{\mY} (\mL f)  g \,d\pi = \int_{\mY} f (\mL^* g) \,d\pi\PERIOD
\]

We next assume that the operator~$-\mL$ is invertible on~$L^2_0(\pi)$, and that it stabilizes some dense space of smooth functions which, together with their derivatives, satisfy some growth conditions (see for instance Theorem~\ref{thm:L-stability} below for Langevin dynamics). When~$\pi$ admits moments of sufficiently high order, it can then be shown (see the discussion in~\cite[Remark~5.5]{lelievre2016partial}) that the linear response~\eqref{eqn:lin-response} is well-defined and can be rewritten, for $f \in L^2(\pi)$ given, as 
\begin{equation}
  \label{eqn:lin-response-2}
  \LIN (f) = -\int_{\mY} \LOPERTIL \mL^{-1} \Pi f \, d\pi \PERIOD
\end{equation}
The above expression serves as the starting point for various reformulations of the linear response that lead to different numerical strategies, as made precise below. These reformulations can be seen as alternatives to the straightforward method which consists in approximating~\eqref{eqn:lin-response} by estimating averages with respect to the invariant measure as time averages over one long realization of the dynamics under consideration and computing, for a given small value of~$\eta$, the estimator
\[
\frac{1}{\eta} \left(\frac1t \int_0^t f(y_s^\eta)\, ds - \frac1t \int_0^t f(y_s) \, ds\right)\,.
\]
This approach corresponds to the so-called Nonequilibrium Molecular Dynamics approach in computational statistical physics~\cite{CKS05,evans-moriss}. When a central limit theorem applies to the time averages under consideration, it is easily seen that the latter estimator has a variance which scales as~$1/(\eta^2 t)$, so that the integration time~$t$ has to be quite large since small values of~$\eta$ are necessary to ensure that the response remains in the linear regime. This indeed motivates turning to alternative expressions of the linear response.

\subsubsection{Green--Kubo method.} One popular reformulation of~\eqref{eqn:lin-response-2} is based on the operator identity
\begin{equation}\label{eqn:L-inverse}
-\mL^{-1} = \int_0^{\infty} \mathrm{e}^{t \mL} \, dt\COMMA
\end{equation}
which can be given a meaning when~$\mathrm{e}^{t \mL}$ is a bounded operator over~$L^2_0(\pi)$, with an operator norm which decays sufficiently fast with respect to~$t$. This allows us to rewrite the linear response as an integrated correlation function, the famous Green--Kubo formula (see for instance~\cite{leimkuhler2016computation, lelievre2016partial}): for $f \in L^2_0(\pi)$,
\begin{equation}
  \label{eqn:GK}
  \LIN (f) = \int_{\mY} \left(-\mL^{-1} f\right)\left(\LOPERTIL^* \mathbf{1}\right) d\pi
  = \int_0^{\infty} \bE^{\pi}\left[ f(y_t) \varphi(y_0) \right] dt,
\end{equation}
where $\varphi = \LOPERTIL^* \mathbf{1}$, and the expectation $\bE^{\pi}$ is taken over the initial conditions distributed as $y_0 \sim \pi$ and over all realizations of the reference dynamics~\eqref{eqn:SDE}. We discuss more precisely the numerical efficiency of the Green-Kubo (GK) estimator based on~\eqref{eqn:GK} in Section~\ref{sec:LR-vs-GK} (in the context of the Langevin dynamics, but our analysis can easily be extended to other diffusion processes). 

\subsubsection{Likelihood ratio method.} A different approach to estimating~\eqref{eqn:lin-response-2} is based on the so-called likelihood ratio formula, obtained by re-weighting the realizations of the reference dynamics with a Girsanov weight to account for the external perturbation, and formally passing first to the small forcing limit, and then to the infinite time limit. As shown in \cite{glynn2019likelihood, plechac2019convergence, wang2019steady}, this leads for $f \in L^2_0(\pi)$ to
\begin{equation}
  \label{eqn:LR}
  \LIN (f) = \lim_{t \to \infty} \bE^{\mu_0}\left[\left(\frac{1}{t}\int_0^t f(y_s) \, ds\right) z_t \right],
  \qquad
  z_t = \int_0^t u(y_s)^{\T} \, dW_s,
\end{equation}
where the function~$u:\mY \to \mathbb{R}^k$ satisfies
\begin{equation}
  \label{eq:su=F}
  \sigma u = F,
\end{equation}
and the expectation $\bE^{\mu_0}$ is taken over initial conditions from an arbitrary initial distribution $\mu_0$, \emph{i.e.}, $y_0 \sim \mu_0$, and over all realizations of the reference dynamics~\eqref{eqn:SDE}. It is clear that~$z_t$ is a zero mean $\mathcal{F}_t$-local martingale.

The main aim of this work is to propose appropriate discrete counterparts of~$z_t$, which we call {\it an auxiliary martingale process}. Another aim is to compare the estimators~\eqref{eqn:GK} and~\eqref{eqn:LR}.

\subsection{Formal strategy for constructing martingale product estimators}
\label{subsec:formal-derivation}

Before proceeding with the rigorous analysis of the numerical methods that we propose for Langevin dynamics, we outline in this section the heuristic derivation of estimators motivated by the likelihood ratio method~\eqref{eqn:LR}. We leave the key question of the construction of the discrete counterparts to the auxiliary martingale process~$z_t$ to Section~\ref{sec:construction_discrete_martingale}. 
However, we warn the reader that these discrete counterparts
are {\it not} based on a direct discretization of the process $z_t$, but on some discrete martingale fluctuation identity -- hence the name {\it martingale product estimators} for the estimation of the linear response. The formal arguments presented in this section are made rigorous for Langevin dynamics in Section~\ref{sec:MP_for_Langevin}.

We consider numerical schemes associated with a fixed timestep~$\dt>0$, written as
\begin{equation}
  \label{eqn:numerstep}
  y^{n+1} = y^n + \Phi_\dt(y^n;G^n)\COMMA
\end{equation}
for a given increment function~$\Phi_\dt$ and random increments~$G^n$ (typically standard Gaussian random variables). We assume that this Markov chain admits a unique invariant probability measure~$\pi_\dt$. Furthermore, the chain is exponentially ergodic provided that the discrete semigroup induced by the Markov chain satisfies both a Lyapunov condition and a minorization condition~\cite{hairer2011yet, meyn2012markov,DMPS18}.  

The basic idea of our approach is to find a zero mean $\mathcal{F}^n$-martingale $\{z^n\}_{n\geq 0}$ re-weighting the centered ergodic average so that
\begin{equation}\label{eqn:desired-estimator}
  \LIN (f) = \lim_{\dt \to 0}\lim_{N \to \infty}\bE_{\dt}\left\{\left(\frac{1}{N}\sum_{n=0}^{N-1} f(y^n) - \pi_{\dt}(f)\right) z^N\right\}\COMMA
\end{equation}
where the expectation $\bE_{\dt}$ is taken over an arbitrary initial distribution~$\mu_0$ for~$y^0$ and over all realizations of the Markov chain~\eqref{eqn:numerstep}. The order of the limits on the right hand side of the above equality is unimportant. 
We comment on the fact that we present the method for the ``perfect'' centering~$\pi_{\dt}(f)$, whereas in practice one should rather center the estimator with the empirical mean of $\{f(y^n)\}_{0 \leq n \leq N-1}$ (see the discussion in~\cite{plechac2019convergence}). 
This issue does not impact the construction of the method. In view of~\eqref{eqn:desired-estimator}, the discovery of a MP estimator boils down to identifying the corresponding discrete auxiliary martingale $\{z^n\}_{n\geq 0}$ and analyzing asymptotic expansions in~$\dt$ and~$1/N$ of the error
\begin{equation}
  \label{eqn:error}
  \bE_{\dt}\left\{\left(\frac{1}{N}\sum_{n=0}^{N-1} f(y^n) - \pi_{\dt}(f)\right) z^N\right\} - \LIN(f)\PERIOD
\end{equation}
We recall below how to quantify time-step biases on the computation of averages along one realization, and then proceed to the analysis of the time-step bias in the product of such averages with discrete martingales as in~\eqref{eqn:error}. Before we do so, we first rewrite the linear response using the solution to the Poisson equation associated with the continuous time process
\begin{equation}
  \label{eqn:cont-Poisson-eqn-general}
  - \mL \widehat{f} = f - \pi(f)\PERIOD
\end{equation}
The linear response~\eqref{eqn:lin-response-2} can then be rewritten as
\begin{equation}
  \label{eq:Lin_f_Poisson}
  \LIN(f) = \int_\mY F^\T \nabla \widehat{f} \, d\pi\,.
\end{equation}

\subsubsection{Quantifying errors on averages along one realization.}
The discrete evolution operator~$P_\dt$ associated with~\eqref{eqn:numerstep} is 
\begin{equation}
  \label{eqn:discsemigroup}
  \left(P_{\dt} f\right)(y) = \bE_{\dt}\{f(y^{n+1}) \, | \,  y^n = y\} = \bE_G\{f(y + \Phi_h(y;G))\}\PERIOD
\end{equation}
Using Taylor expansions, the action of the operator~$P_h$ can usually be written as an expansion in powers of~$\dt$: for any $C^{\infty}$ test function $f$,
\begin{equation}
  \label{eq:expansion_Pdt}
  P_{\dt}f = f + \dt \mA_1 f + \ldots + \dt^{\alpha} \TAOP_{\alpha}f + \dt^{\alpha+1} \mathcal{R}_{\alpha, \dt}f \COMMA
\end{equation}
where~$\TAOP_i$ and~$\mathcal{R}_{\alpha, \dt}$ are operators which depend on the underlying numerical scheme $\Phi_h$ (the operators~$\TAOP_i$ being linear differential operators with smooth coefficients). Since the semigroup of the continuous time process~\eqref{eqn:SDE} has the expansion
\[
    \mathrm{e}^{\dt \mL} f = f + \dt \mL f + \frac{\dt^2}{2} \mL^2 f + \dots +
    \frac{\dt^{\alpha}}{\alpha !}\mL^\alpha f + \BIGO(\dt^{\alpha+1}) \COMMA
\]
it is easily seen that $\mA_1 = \mL$ corresponds to weak first order schemes ($\alpha=1$), while the additional condition $\mA_2 = \mL^2/2$ characterizes weak second order schemes $(\alpha=2$); see, for instance, \cite{KP92,MT04} for an introduction to the numerical analysis of discretization schemes for diffusion processes.
Moreover, it can be shown that an ergodic scheme of weak order~$\alpha$ has an invariant probability measure which is correct at least at order~$\alpha$, in the following sense: for any smooth function~$f$ satisfying some growth conditions, there exists $C_f \in \mathbb{R}_+$ such that 
\begin{equation}
  \label{eq:error_inv_meas_dt}
  \left|\pi_\dt(f) - \pi(f)\right| \leq C_f \dt^{\alpha}.
\end{equation}

We can now recall the approximation result for pathwise averages of~$\pi(f)$ using a numerical scheme of weak order~$\alpha$ (see for instance~\cite{mattingly2010convergence}):
\begin{equation}
  \label{eqn:ergodicaverage}
  \left|\bE_{\dt}\left\{ \frac{1}{N}\sum_{n=0}^{N-1} f(y^n) \right\} - \pi(f)\right| \leq C\left(\dt^{\alpha} + \frac{1}{N\dt}\right)\PERIOD
\end{equation}
The proof of this result is based on a decomposition of the pathwise average into a sum of discrete martingale increments. To make this decomposition precise, since it will be of paramount importance for the numerical analysis we present, we introduce the operator 
\begin{equation}
  \label{eqn:discrete_generator}
  \LOPERH = \frac{1}{\dt}(P_\dt - I)\PERIOD
\end{equation}
By definition, 
\begin{equation}
  \label{eqn:obs_increment}
  f(y^{n+1}) - f(y^n) = \dt \LOPERH(f)(y^n) + \DMN(f)\COMMA
  \qquad
  \DMN(f) = f(y^{n+1}) - \left(P_\dt f\right)(y^n)\COMMA
\end{equation}
with $\DMN(f)$ the increment (\emph{i.e.}, martingale difference) of a discrete martingale. We next define the discrete Poisson equation
\begin{equation}
  \label{eqn:disc-Poisson-eqn-general}
  -\LOPERH \WHF = \frac{1}{\dt}(I - P_{\dt}) \WHF = f - \pi_{\dt}(f)\PERIOD
\end{equation}
This equation is well-posed when~$P_\dt^n$ decays sufficiently fast with~$n$ in a suitable functional space, typically spaces of measurable functions with a maximal growth at infinity dictated by a Lyapunov function, in which case (see for instance~\cite{meyn2012markov})
\[
(I - P_{\dt})^{-1} = \sum_{n=0}^{\infty} P_{\dt}^n \PERIOD
\]
Using the solution~$\WHF$ of~\eqref{eqn:disc-Poisson-eqn-general}, we can rewrite pathwise averages as
\begin{equation}\label{eqn:estimator-idea1}
\begin{split}
\frac{1}{N}\sum_{n=0}^{N-1} f(y^n) - \pi_{\dt}(f) & = 
-\frac{1}{N}\sum_{n=0}^{N-1} \LOPERH \WHF(y^n) 
 = \frac{1}{N\dt} \sum_{n=0}^{N-1} \left[\DMN(\WHF)+\WHF(y^n)-\WHF(y^{n+1})\right]\\
&= \frac{1}{Nh} \left[ \WHF\left(y^0\right) - \WHF\left(y^N\right)\right] + \frac{1}{Nh}\sum_{n=0}^{N-1} \DMNFH   \PERIOD
\end{split}
\end{equation}
The error estimate~\eqref{eqn:ergodicaverage} then follows by taking the expectation of both sides of the previous equality (provided~$\bE_{\dt}\{\WHF(y^n)\}$ can be bounded uniformly in~$n \geq 0$), and recalling~\eqref{eq:error_inv_meas_dt}. 

\subsubsection{Quantifying errors in~\eqref{eqn:error}.}
Following the derivation leading to~\eqref{eqn:estimator-idea1}, we start by rewriting the first term of~\eqref{eqn:error} as
\begin{equation}
  \label{eqn:estimator-idea}
  \begin{split}
    \bE_{\dt}\left\{\left(\frac{1}{N}\sum_{n=0}^{N-1} f(y^n) - \pi_{\dt}(f) \right) z^N\right\} 
    & =
    \frac{1}{N\dt}\bE_{\dt}\left\{\sum_{n=0}^{N-1} \DMNFH    z^N\right\} \\
    & \ \ + 
    \frac{1}{N\dt}\bE_{\dt}\left\{  \left(\WHF(y^0) - \WHF(y^N)\right) z^N\right\}\,.
  \end{split}
\end{equation}
At this point, it is tempting to expand $\DMN(\WHF)$ in powers of~$\dt$ using a Taylor expansion of~$\WHF(y^{n+1})=\WHF(y^{n}+\Phi_\dt(y^n;G^n))$ at $y^n$. However, such an expansion requires sufficient regularity of the solution~$\widehat{f}_\dt$ to the discrete Poisson equation~\eqref{eqn:disc-Poisson-eqn-general}. A control of even the first derivatives of the solution to the discrete Poisson equation is not guaranteed in general since the natural functional spaces for the well posedness of the discrete Poisson equation are (weighted) spaces of bounded measurable functions. However, there exists a sufficiently regular function~$\WTF$ which approximates the solution~$\widehat{f}_\dt$ to the discrete Poisson equation to an arbitrary order in powers of~$\dt$; and also happens to approximate the solution to the continuous time Poisson equation~\eqref{eqn:cont-Poisson-eqn-general}. Indeed, it turns out that, by a technical result whose proof is deferred to Section~\ref{subsec:approx-Poisson-equation} for Langevin dynamics (see also~\cite{plechac2019convergence} for non-degenerate diffusions on a torus), we obtain the following: there exists a function~$\WTF$, and operators~$\mathcal{Q}_1$ and~$\mathcal{Q}_2$ such that
\begin{equation}
  \label{eqn:approx-operator-preview}
    \WTF - \widehat{f} 
    =
    \left\{
    \begin{array}{ll}
    -\dt \mathcal{Q}_1(f - \pi(f)),     &  \alpha = 1;\\[3pt]
    \dfrac{\dt}{2}(f - \pi(f)) -\dt^2 \mathcal{Q}_2(f - \pi(f)),   &  \alpha =2\COMMA
    \end{array}
    \right.
\end{equation}
where $\alpha=1$ or~$2$ refers to the order of the weak approximation order for the numerical discretization~\eqref{eqn:numerstep} at hand.

We next replace $\widehat{f}_\dt$ by~$\WTF$ in the sum on the right hand side of~\eqref{eqn:estimator-idea} (with a remainder which can be as small as needed, but we restrict ourselves to~$\BIGO(\dt^{\alpha})$ since this will be sufficient for our analysis), thus obtaining
\begin{equation}
  \label{eqn:general-MP-estimate}
  \begin{aligned}
    \bE_{\dt}\left\{\left(\frac{1}{N}\sum_{n=0}^{N-1} f(y^n) - \pi_{\dt}(f) \right) z^N\right\} 
    & =
    \frac{1}{N\dt}\bE_{\dt}\left\{\sum_{n=0}^{N-1} \DMNFT  z^N\right\} \\
    & \quad + 
    \frac{1}{N\dt}\bE_{\dt}\left\{  \left(\widehat{f}_\dt(y^0) - \widehat{f}_\dt(y^N)\right) z^N\right\}
    + 
    \BIGO(\dt^{\alpha})\PERIOD
  \end{aligned}
\end{equation}
Ensuring that the remainder term in the previous estimate is indeed of order~$\dt^\alpha$ requires some work; see for instance~\eqref{eqn:1st-order-estimate}-\eqref{eqn:1st-order-estimate-1}.
Under mild conditions, the second term on the right hand side of~\eqref{eqn:general-MP-estimate} can be controlled by~$\BIGO((N\dt)^{-1/2})$ uniformly with respect to the timestep~$\dt$. 
This term therefore contributes to the finite integration time error (it is larger than for standard time averages, compare with the term~$\BIGO((N\dt)^{-1})$ in~\eqref{eqn:ergodicaverage}), and hence the dominating discretization error is completely determined by the first term on the right hand side of~\eqref{eqn:general-MP-estimate}. In view of~\eqref{eq:Lin_f_Poisson}, our goal is therefore to construct the auxiliary discrete martingale~$z^N$ by exploring the expansion in~$\dt$ of the error term
\[
\frac{1}{N\dt}\bE_{\dt}\left\{\sum_{n=0}^{N-1} \DMNFT  z^N\right\} - \int_\mY F^\T \nabla \widehat{f} \,d\pi\PERIOD
\]

\subsection{Construction of the discrete auxiliary martingale}
\label{sec:construction_discrete_martingale}

As mentioned above, rather than starting with the continuous auxiliary martingale~$\{z_t\}_{t\geq 0}$, we construct a discrete time auxiliary martingale~$\{z^n\}_{n\geq 0}$ which allows to obtain a consistent approximation of~$\LIN(f)$ in the limit~$\dt\to 0$. We consider a specific form for the auxiliary martingale, namely 
\begin{equation}
  \label{eq:def_z_n}
  z^{n+1} - z^n = \DMN(g) \equiv g(y^{n+1}) - \left(P_{\dt} g\right)(y^n), \qquad z^0 = 0\COMMA
\end{equation}
where~$g$ is a sufficiently regular function chosen so that
\begin{equation}
  \label{eq:estimator_Lin_f}
\LIN_{\dt,N}(f) \equiv \frac{1}{N\dt}\bE_{\dt}\left\{\sum_{n=0}^{N-1} \DMNFT  \sum_{m=0}^{N-1} \DMM(g)\right\}
\end{equation}
approximates $\LIN(f)$ to the desired order~$\BIGO(\dt^\alpha)$. Clearly, $\{z^n\}_{n\geq 0}$ defined in~\eqref{eq:def_z_n} is a zero mean martingale adapted to~$\mathcal{F}^n$. In fact, since $\DMN(\WTF)$ and~$\DMM(g)$ are increments of mean zero which are conditionally independent for $m\ne n$, 
\[
\frac{1}{N\dt}\bE_{\dt}\left\{\sum_{n=0}^{N-1} \DMNFT  \sum_{m=0}^{N-1} \DMM(g)\right\} = \frac{1}{N\dt} \sum_{n=0}^{N-1} \bE_{\dt}\left\{\DMNFT \DMN(g)\right\}.
\]
The expression on the right hand side of the previous equality motivates the name ``martingale product'' of our estimator.

The next step is to obtain an expansion in powers of~$\dt$ of the right hand side of the previous equality. To simplify the notation and the computation of these expansions, we introduce the ``carr\'e du champ operator'' associated with the generator~$\mL$: for two smooth funtions~$f,g$, 
\begin{equation}
  \label{eqn:GammaDef}
  \Gamma(f,g) = \mL(fg) - f\mL(g) - g\mL(f) = \nabla f^\T S \nabla g\PERIOD
\end{equation}
We also consider its discrete counterpart (recall the definition~\eqref{eqn:discrete_generator} for~$\LOPERH$):
\begin{equation}
  \label{eqn:GammahDef}
  \Gamma_\dt(f,g) = \LOPERH(fg) - f\LOPERH(g) - g\LOPERH(f) - \dt\LOPERH(f)\LOPERH(g)\PERIOD
\end{equation}
In view of the definition~\eqref{eqn:GammahDef}, the following crucial identity holds
\begin{equation}
  \label{eq:crucial_identity}
  \bE_{G^n}\left\{\DMNFT \DMN(g)\right\}
  = \left[P_{\dt}(\WTF g)\right](y^n) - (P_{\dt}\WTF)(y^n) (P_{\dt}g)(y^n)
  = \dt \Gamma_\dt(\WTF,g)(y^n)\COMMA
\end{equation}
where $\bE_{G^n}$ denotes the expectation with respect to the random variable~$G^n$ in~\eqref{eqn:numerstep}. Therefore, for a numerical scheme of weak order~$\alpha$,
\begin{equation}
  \label{eqn:expansionDMDN}
  \LIN_{\dt,N}(f) = \frac{1}{N} \sum_{n=0}^{N-1}\bE_\dt\left\{\Gamma_h(\WTF,g)(y^n)\right\} = \int_\mY \Gamma_h(\WTF,g) \, d\pi + \BIGO(\dt^\alpha) + \BIGO\left(\frac{1}{N\dt}\right)\COMMA
\end{equation}
where the second equality follows from applying~\eqref{eqn:ergodicaverage} to the observable $\Gamma_h(\WTF, g)$. In order to expand $\LIN_{\dt,N}(f)$ in powers of~$\dt$, it suffices to expand the integral in the last equality of~\eqref{eqn:expansionDMDN}. We consider successively schemes of weak order~$\alpha = 1$ or~2, using the following expansions:
\begin{eqnarray}
\mbox{$\alpha=1$:} & \LOPERH = \mL + \BIGO(h)\COMMA
&\Gamma_\dt(f,g) = \Gamma(f,g) + \BIGO(h)\COMMA \label{eqn:order1}\\
\mbox{$\alpha=2$:}  & \dps \LOPERH = \mL + \frac{h}{2} \mL^2 + \BIGO(h^2)\COMMA
 &\Gamma_\dt(f,g) = \Gamma(f,g) + h \Lambda(f,g) + \BIGO(h^2)\COMMA \label{eqn:order2}
\end{eqnarray}
where
\begin{equation}
  \label{eqn:Lambda-operator}
\Lambda(f, g) = \frac{1}{2}\Big(\mL^2(fg) - f \mL^2 g - g\mL^2 f -  2 (\mL f)( \mL g)\Big)\PERIOD
\end{equation}
Note that
\begin{equation}
  \label{eq:int_Lambda_pi}
  \begin{aligned}
  \int_\mY \Lambda(f, g) \, d\pi & = -\int_\mY \Big(\frac{\mL+\mL^*}{2}\Big)f \mL g + \mL f \Big(\frac{\mL+\mL^*}{2}\Big) g\, d\pi \\
  & = \frac12 \int_\mY \nabla f^\T S \nabla(\mL g) + \nabla g^\T S \nabla(\mL f) \, d\pi\PERIOD
  \end{aligned}
\end{equation}

\subsubsection{Weak first order schemes ($\alpha=1$).}
\label{sec:first_order_general}

For a weak first order discretization scheme, \eqref{eqn:approx-operator-preview}, \eqref{eqn:GammaDef}, \eqref{eqn:expansionDMDN} and~\eqref{eqn:order1} lead to
\begin{equation}
  \label{eqn:carre-du-champ-1st-order}
  \LIN_{\dt,N}(f) = \int_\mY \Gamma(\widehat{f}, g) \, d\pi + \BIGO(\dt) + \BIGO\left(\frac{1}{N\dt}\right) = \int_\mY  \nabla \widehat{f}^\T S \nabla g \, d\pi + \BIGO(\dt) + \BIGO\left(\frac{1}{N\dt}\right)\PERIOD
\end{equation}
This coincides at dominant order with the linear response as given by~\eqref{eq:Lin_f_Poisson} provided the following condition holds: for all smooth functions~$\varphi$ with compact support,
\begin{equation}
  \label{eqn:constraint-g}
  \int_\mY \nabla \varphi^\T S \nabla g \, d\pi = \int_\mY \nabla \varphi^\T F \, d\pi\,.
\end{equation}
The latter condition can be seen as the Euler--Lagrange equation associated with the minimization problem
\[
\min_{\phi \in H^1(\pi)} \left\{ \int_\mY \left|\sigma^\T\nabla \phi-u\right|^2 d\pi \right\},
\]
where~$u$ satisfies~\eqref{eq:su=F}, \emph{i.e.}, $\sigma u = F$, and $H^1(\pi) = \{\phi\in L^2(\pi) ~|~ \partial_{y_i} \phi \in L^2(\pi), \forall i = 1, \ldots, d\} $. This minimization problem always admits a solution when~$\sigma$ is bounded (since it can be seen as the orthogonal projection of~$u$ onto the closed subspace~$\{\sigma^\T \nabla \phi ~| \ \phi \in H^1(\pi)\}$ of~$L^2(\pi)$), but the corresponding minimizer~$g$ may not be unique, even if the function~$g$ is required to have mean zero with respect to~$\pi$ (in order to remove the trivial obstruction to uniqueness provided by adding constants to a minimizer). Alternatively, \eqref{eqn:constraint-g} can be interpreted as a weak formulation of the following Poisson-type equation for~$g$:
\begin{equation}
  \label{eq:PDE_g}
  \nabla^* (S \nabla g) = \nabla^* F \equiv \sum_{i=1}^d \partial_{y_i}^* F_i,
\end{equation}
with $\partial_{y_i}^* = -\partial_{y_i} - \partial_{y_i} (\ln \rho)$ the adjoint of~$\partial_{y_i}$ with respect to the canonical scalar product on~$L^2(\pi)$.
The well-posedness of the equation defining~$g$ is discussed in Section~\ref{sec:LD-application} for Langevin dynamics. When~\eqref{eqn:constraint-g} holds, we can then consider (at least formally) the following first order MP estimator
\begin{equation}
  \label{eqn:general-first-order-MP}
  \mathcal{E}_{\dt, N}^{\mathrm{MP1}} = \frac{1}{N}\sum_{n=0}^{N-1}\left( f(y^n) -\pi_{\dt}(f)\right) \left( g(y^{n+1}) - (P_{\dt} g)(y^n) \right)\PERIOD
\end{equation}
To obtain an actual expression of~$z^n$ which can be implemented in practice, the final step is to Taylor expand $g(y^{n+1})$ in (fractional) powers of~$\dt$ around~$y^n$, and to replace~$P_\dt$ in~$(P_{\dt} g)(y^n)$ with~\eqref{eq:expansion_Pdt}.

\begin{example}[Euler-Maruyama scheme]
{\rm 
Consider a non-singular noise ($\sigma(y)$ is invertible for all~$y \in \mY$), and the Euler-Maruyama discretization 
\[
y^{n+1} = y^n + h\, b(y^n) + h^{1/2}\, \sigma(y^n) G^n\COMMA 
\]
where $G^n$ are independent and identically distributed standard random Gaussian variables with identity covariance. A simple computation then shows that the martingale increment appearing in the definition~\eqref{eq:def_z_n} can be expanded as~$\DMN(g) = h^{1/2} \nabla g(y^n)^\T\sigma(y^n) G^n +\BIGOM(\dt^{3/2})$, 
~where $\BIGOM(\dt^\alpha)$ denotes random terms of order~$h^\alpha$ (in the sense that $\bE\{|\BIGOM(\dt^\alpha)|\} = \BIGO(\dt^\alpha)$). The condition~\eqref{eqn:constraint-g} suggests to replace~$\sigma^\T\nabla g$ with~$u=\sigma^{-1}F$, which leads to the following definition for the discrete auxiliary martingale: 
\[
z^{n+1} = z^n + h^{1/2} \left(G^n\right)^\T \sigma^{-1}(y^n) F(y^n)\PERIOD
\]
The discrete martingale~$z^n$ is in this case a consistent approximation of the likelihood process~$z_t$, and the associated first order estimator is 
\[
\mathcal{E}_{\dt, N}^{\mathrm{MP1}} = \frac{1}{N}\sum_{n=0}^{N-1}\left( f(y^n) -\pi_{\dt}(f)\right) \left( \dt^{1/2} \left(G^n\right)^\T \sigma^{-1}(y^n) F(y^n)  \right)\PERIOD
\]
This corresponds to the first order estimator already introduced in~\cite{plechac2019convergence}.
}
\end{example}

\subsubsection{Weak second order schemes $\alpha=2$.}
The interest of our approach becomes more apparent with a second order discretization scheme as it allows us to systematically discover the corresponding second order MP estimator. In order to correct the discretization error of order~$\dt$, we consider an additional martingale of order~$\BIGOM(\dt^{3/2})$, which can be seen as a correction to~\eqref{eq:def_z_n}, so that
\[
z^{n+1} - z^n = \DMN(g+\dt\WTG ) \COMMA
\]
where the leading order function~$g$ is determined by~\eqref{eqn:constraint-g} as for first order schemes.
From~\eqref{eqn:approx-operator-preview}, \eqref{eqn:expansionDMDN},~\eqref{eqn:order2} and the bilinearity of the operator $\Gamma$, 
\[
  \begin{aligned}
    \LIN_{\dt,N}(f) & = \int_\mY \Gamma_\dt\left(\widehat{f}+\frac{\dt}{2}f, g+\dt\WTG \right) d\pi + \BIGO(\dt^2) + \BIGO\left(\frac{1}{N\dt}\right) \\
    & = \int_\mY \Gamma\left(\widehat{f},g\right) d\pi + \dt \int_\mY \frac12 \Gamma(f,g) + \Gamma\left(\widehat{f},\WTG \right) + \Lambda\left(\widehat{f},g\right) d\pi+ \BIGO(\dt^2) + \BIGO\left(\frac{1}{N\dt}\right)\PERIOD
  \end{aligned}
\]
As for first order schemes, this coincides at leading order with the linear response as given by~\eqref{eq:Lin_f_Poisson} provided~\eqref{eqn:constraint-g} holds. The bias of order~$\dt$ vanishes provided~$\WTG $ is chosen such that
\[
  \begin{aligned}
    \int_\mY \Gamma\left(\widehat{f},\WTG \right) d\pi & = - \int_\mY \frac12 \Gamma(f,g) + \Lambda\left(\widehat{f},g\right) d\pi \\
    & = - \frac12 \int_\mY \nabla f^\T S \nabla g \, d\pi - \frac12 \int_\mY \nabla (\mL \widehat{f})^\T S \nabla g \, d\pi - \frac12\int_\mY \nabla \widehat{f}^\T S \nabla (\mL g) \, d\pi \\
    & = - \frac12 \int_\mY \nabla \widehat{f}^\T S \nabla (\mL g) \, d\pi = -\frac12\int_\mY \Gamma(\widehat{f},\mL g) \, d\pi \COMMA
  \end{aligned}
\]
where we used~\eqref{eq:int_Lambda_pi} to obtain the second equality and ~\eqref{eqn:cont-Poisson-eqn-general} to obtain the third equality. This calculation suggests to choose
\[
\WTG  = -\frac12 \mL g\,.
\]
This leads to the following second order MP estimator:
\begin{equation}
  \label{eqn:general-second-order-MP}
  \mathcal{E}_{\dt, N}^{\mathrm{MP}2} = \frac{1}{N}\sum_{n=0}^{N-1}\left[ f(y^n) -\pi_{\dt}(f)\right] \left[ \left(g-\frac{\dt}{2}\mL g\right)(y^{n+1}) - P_{\dt} \left(g-\frac{\dt}{2}\mL g\right)(y^n) \right]\COMMA
\end{equation}
where $g$ satisfies~\eqref{eqn:constraint-g}. As discussed after~\eqref{eqn:general-first-order-MP}, one needs to Taylor expand $g(y^{n+1})$ in (fractional) powers of~$\dt$ around~$y^n$, and to replace~$P_\dt$ in~$(P_{\dt} g)(y^n)$ with~\eqref{eq:expansion_Pdt} in order to obtain a computable expression for $z^n$. This answers in any case the question left open in~\cite{plechac2019convergence} about constructing second order schemes for general non degenerate diffusion processes in spaces of dimension higher than or equal to two.

We emphasize that the derivation of the MP estimators~\eqref{eqn:general-first-order-MP} and~\eqref{eqn:general-second-order-MP} presented in this section is 
heuristic since we have not provided any uniform controls for the error terms. In order to carry out rigorous analysis in the next sections
we consider the particular case of Langevin dynamics with numerical schemes based on operator splitting methods~\cite{leimkuhler2013rational}.

\section{Splitting schemes for Langevin dynamics}
\label{sec:LD-application}

In this section, we present Langevin dynamics and define the class of perturbations we consider (see Section~\ref{sec:Langevin_transport}). We next recall in Section~\ref{sec:splitting_Lang} splitting schemes to discretize Langevin dynamics, as well as some of their properties. We conclude in Section~\ref{subsec:approx-Poisson-equation} by providing a construction of the smooth approximations to discrete Poisson equations introduced in~\eqref{eqn:approx-operator-preview}.

\subsection{Langevin dynamics and transport coefficients}
\label{sec:Langevin_transport}

Langevin dynamics is a diffusion process which describes the evolution of $\cN$ particles in a space of physical dimension~$d$, at fixed temperature. We denote by $q = (q_{1}, \cdots, q_{\cN}) \in \mathcal{D}$ the positions of the particles (typically, $\mathcal{D} = \mathbb{T}^D$ with $\mathbb{T} = \mathbb{R}\backslash \mathbb{Z}$ and~$D=\cN d$, or~$\mathcal{D}=\mathbb{R}^D$), and by~$p = (p_{1}, \cdots, p_{\cN})\in \mathbb{R}^{D}$ their momenta. The phase space for the elements~$y=(q,p)$ is therefore $\mY = \mathcal{D} \times  \mathbb{R}^{D}$.

\subsubsection{Reference Langevin dynamics.}
Langevin dynamics is defined by the system of stochastic differential equations
\begin{equation}
  \label{eqn:ULD}
  \left\{
  \begin{aligned}
      dq_t &= M^{-1}p_t\,dt\,, \\
      dp_t &= -\nabla V(q_t)\, dt - \gamma M^{-1} p_t \,dt + \displaystyle\sqrt{\frac{2\gamma}{\beta}} dW_t\,,
  \end{aligned} \right.
\end{equation}
where~$\beta$ is (proportional to) the inverse temperature, $V:\mathcal{D} \to \mathbb{R}$ is the potential function, the mass matrix $M$ is a symmetric definite positive matrix such as $M= \mathrm{diag}(m_1 \Id_d, \ldots, m_N \Id_d)$ with~$m_i > 0$ for $1 \leq i \leq \cN$, $\beta$ is proportional to the inverse temperature, $\gamma>0$ is the friction coefficient, and $W_t$ is a standard $D$-dimensional Brownian motion. We consider the following running assumption.

\begin{assumption}
\label{ass:V_smooth_and_domain_compact}
  The potential $V$ belongs to~$C^{\infty}(\mathcal{D})$ and the position space is compact: $\mathcal{D} = \mathbb{T}^D$.
\end{assumption}

Periodic boundary conditions for the positions are natural in molecular dynamics simulations of condensed matter systems~\cite{FrenkelSmit,tuckerman2010statistical}, so that the assumption $\mathcal{D} = \mathbb{T}^D$ is physically relevant. The main benefit of this assumption is that it allows us to simplify various mathematical arguments, although many of them could be extended to unbounded position spaces.

The existence and uniqueness of strong solutions to~\eqref{eqn:ULD} is standard when the position space is compact~\cite{bellet2006ergodic}. Moreover, the unique invariant measure of~\eqref{eqn:ULD} is the Boltzmann--Gibbs probability measure
\[
\pi(dq\, dp) = Z^{-1} \mathrm{e}^{-\beta H(q, p)}\, dq\,dp\,,\; \qquad Z = \int_\mY \mathrm{e}^{-\beta H}\,, \;\;\qquad H(q, p) = V(q) + \frac{1}{2}p^{\T} M^{-1} p\,,
\] 
where $H$ denotes the Hamiltonian of the system.

\subsubsection{Poisson equation associated with the reference Langevin dynamics.}
In order to control remainder terms we restrict the analysis to a subclass class of smooth observables.
Since the position space~$\mathcal{D}$ is compact, it is sufficient to control the growth of functions and their derivatives with respect to the momentum variable. Specifically, as in~\cite{talay2002stochastic,kopec2015weak}, we consider scale functions of polynomial form
\[
  \mKS(q, p) = 1 + |p|^{2s}\,, \;\;\;\qquad s = 1, 2, \ldots\COMMA
\]
and define
\[
B_{\mKS}^{\infty}(\mY) = \left\{ f \mathrm{~measurable} \ \middle| \ \left\|f\right\|_{B_{\mKS}^{\infty}} = \sup_{y \in \mY} \left|\frac{f(y)}{\mKS(y)}\right| < \infty \right\}\,.
\]
The set of observables we consider is then composed of smooth functions which, together with their derivatives, grow at most polynomially:
\[
\mS = \left\{f \in C^{\infty}(\mY) \left| ~\forall k \in \mathbb{N}^{2D}, \ \exists s \in \mathbb{N} ~\text{such that}~ \partial^k f \in B_{\mKS}^{\infty}(\mY) \right. \right\}\COMMA
\]
where $\partial^k f = \partial_{y_1}^{k_1}\dots \partial_{y_{2D}}^{k_{2D}} f$ for $k = (k_1,\dots,k_{2D}) \in \mathbb{N}^{2D}$. Clearly, the invariant measure $\pi$ integrates all functions~$\mKS$ for $s \in \mathbb{N}$, and hence $\mS$ is a dense subspace of~$L^2(\pi)$. Finally, we define
\[
\mS_0 = \Pi \mS = \mS \cap L^2_0(\pi)\,.
\]

An important result for our analysis is the following regularity result concerning the existence and uniqueness of solutions to the continuous time Poisson equation
\begin{equation}
  \label{eqn:cont-Poisson-eqn}
  -\mL_{\gamma} \widehat{f} = f - \pi(f)\COMMA
\end{equation}
where $\mathcal{L}_{\gamma}$ is the generator of the Langevin dynamics: 
\begin{equation}
\mathcal{L}_{\gamma} = A + B + \gamma C \COMMA
\label{eqn:Lgamma}
\end{equation}
with
\begin{equation}
A = p^\T M^{-1} \nabla_q\,,
  \qquad
  B = -\nabla V(q)^{\T} \nabla_p\,,
  \qquad 
  C = - p^\T M^{-1} \nabla_p + \frac{1}{\beta} \Delta_p\,. \label{eqn:ABC}
\end{equation}
The invertibility of~$\mL_{\gamma}$ considered as an operator on~$L_0^2(\pi)$ can be obtained by results of hypocoercivity (relying on the exponential convergence of the semigroup as provided by~\cite{Herau06,villani2009hypocoercivity,DMS15} for instance, see the introduction of~\cite{bernard2020hypocoercivity} for an
extensive review) or by a direct analysis based on Schur complements~\cite{bernard2020hypocoercivity}. In fact, the solution to~\eqref{eqn:cont-Poisson-eqn} belongs in fact to~$\mS_0$ when~$f \in \mS$, as made precise in the following result proved in~\cite{kopec2015weak, talay2002stochastic}.

\begin{theorem}
  \label{thm:L-stability}
  The space $\mS_0$ is stable under $\mL_{\gamma}^{-1}$: %
  for all $f \in \mS$, there is a unique $\widehat{f} \in \mS_0$ %
  satisfying the Poisson equation~\eqref{eqn:cont-Poisson-eqn} associated with $\mL_{\gamma}$.  %
\end{theorem}

The regularity result of Theorem~\ref{thm:L-stability} is of fundamental importance for justifying the consistency of MP estimators, since it allows 
us to replace solutions to the Poisson equation associated with discretizations of the stochastic differential equation~\eqref{eqn:ULD} by smooth functions, up to a small error term (see Section~\ref{subsec:approx-Poisson-equation} below). 

\subsubsection{Perturbations of the Langevin dynamics.}
Various dynamical properties of Langevin dynamics can be of interest, in particular transport coefficients such as the mobility, thermal conductivity or shear viscosity; see for instance~\cite{evans-moriss,tuckerman2010statistical} as well as~\cite{hairer2008ballistic,Lefevere,JS12,BK21} for some representative mathematical studies dealing with these transport coefficients. 
As discussed in Section~\ref{sec:LR_stoch_dyn}, we view transport coefficients as being obtained from the linear response of steady state averages of Langevin dynamics perturbed by an external forcing. More precisely, we consider a smooth external forcing $F: \mathcal{D} \to \mathbb{R}^{D}$, and the following non-equilibrium Langevin dynamics for~$\eta \in \mathbb{R}$:
\begin{equation}
  \label{eqn:ULD-perturbed}
  \left\{
    \begin{aligned}
      dq_t^\eta &= M^{-1}p_t^\eta\,dt \COMMA\\
      dp_t^\eta &= \left(-\nabla V(q_t^\eta) + \eta F(q_t^\eta)\right)\, dt - \gamma M^{-1} p_t^\eta \,dt + \displaystyle\sqrt{\frac{2\gamma}{\beta}} dW_t \PERIOD
    \end{aligned} \right.
\end{equation}
Note that the forcing~$F(q)$ is a function of the positions only, and is in general not the gradient of a smooth periodic function. 

For instance, mobility can be computed with~$F(q)$ constant (as made precise in~\cite[Section~5]{lelievre2016partial}), while the shear viscosity can be recovered by considering a force whose magnitude in one direction depends only on the component of the position in another direction (the so-called sinusoidal transverse field method~\cite{GMS73,JS12}). It can be shown that~\eqref{eqn:ULD-perturbed} admits a unique strong solution. Moreover, this SDE has a unique invariant probability measure~$\pi^{\eta}$, using standard arguments (a minorization condition is obtained as in~\cite{MSH02} and we can then use the results of~\cite{hairer2011yet} since a Lyapunov condition holds as well). In addition, the invariant measure has a smooth density with respect to the Lebesgue measure by hypoellipticity~\cite{hormander1967hypoelliptic}.

\begin{remark}
{
  It would be possible to consider an external forcing which depends both the positions~$q$ and momenta~$p$. In order to guarantee the existence and uniqueness of strong solutions to~\eqref{eqn:ULD-perturbed}, at least for~$|\eta|$ small enough, this would require conditions on~$F$ such that~$\mKS$ are Lyapunov functions. This is for instance the case when there exist $a$, $b \in \mathbb{R}$ such that
  \[
    p^\T F(q,p) \leq a |p|^2 + b\,.
  \]
  However, we refrain from considering such a more general setting since the applications we have in mind consider a position-dependent forcing only.
  }
\end{remark}

The generator of the perturbed Langevin dynamics~\eqref{eqn:ULD-perturbed} can be written as \[
\mL_{\gamma}^{\eta} = \mL_{\gamma} + \eta\LOPERTIL\,,
\]
with $\mathcal{L}_{\gamma}$ defined in~\eqref{eqn:Lgamma}, and the perturbation operator
\begin{equation}
\LOPERTIL = F^\T \nabla_p \PERIOD
\label{eqn:Ltilde}
\end{equation}
Following the argument in~\cite[Remark~5.5]{lelievre2016partial}, it can be shown that the linear response of an observable~$f \in \mS$ is well defined and can be written as
\begin{equation*}
\LIN(f) = -\int_{\mY} \LOPERTIL \mL_{\gamma}^{-1} \Pi f \,d\pi \PERIOD
\end{equation*}
In view of the Poisson equation~\eqref{eqn:cont-Poisson-eqn} associated with the generator~\eqref{eqn:Lgamma}, we end up with the reformulation
\begin{equation}
\label{eqn:Langevin-lin-res-reformulation}
\LIN(f) = \int_{\mY} F^\T \nabla_p \widehat{f}\, d\pi \COMMA
\end{equation}
that serves as the starting point for deriving MP estimators.

\subsection{Splitting schemes for Langevin dynamics and their properties}
\label{sec:splitting_Lang}

We briefly describe in this section splitting schemes to discretize the Langevin dynamics~\eqref{eqn:ULD} (see~\cite{leimkuhler2013rational, leimkuhler2016computation} for a full exposition) and then recall some of their properties, in particular concerning their invariant measures and their convergence to stationarity.

\subsubsection{Splitting schemes.}
The splitting schemes we consider for Langevin dynamics are based on the decomposition \eqref{eqn:Lgamma} of the generator $\mathcal{L}_{\gamma}$.
First order splitting schemes with a timestep~$\dt$ are obtained by a Lie-Trotter approximation of the continuous time evolution operator~$\mathrm{e}^{\dt\mL_{\gamma}}$ by the discrete time evolution operator
\[
P_{\dt}^{Z, Y, X} = \mathrm{e}^{\dt Z} \mathrm{e}^{\dt Y}  \mathrm{e}^{\dt X},
\] 
where $(Z, Y, X)$ is one of the six possible permutations of the elementary operators~$(A, B, \gamma C)$ introduced in~\eqref{eqn:ABC}. For example, the scheme associated with~$P_{\dt}^{B, A, \gamma C}$ reads
\begin{equation}
  \label{eqn:BAC-scheme}
  \left\{
    \begin{aligned}
      p^{n+\frac{1}{2}} &= p^n - \dt \nabla V(q^n) \COMMA \\
      q^{n+1} &= q^n + \dt M^{-1} p^{n + \frac{1}{2}} \COMMA \\
      p^{n+1} &= \alpha_{\dt} p^{n+\frac{1}{2}} + \sqrt{\frac{1-\alpha_{\dt}^2}{\beta} M} G^n,
    \end{aligned}
  \right.
\end{equation}
where $\alpha_{\dt} = \exp(-\gamma M^{-1} \dt)$ and $(G^n)_{n \geq 0}$ are independent and identically distributed standard $D$-dimensional Gaussian random vectors.

Similarly, second order schemes approximate the continuous time evolution operator~$\mathrm{e}^{\dt\mL_{\gamma}}$ through a Strang splitting
\[
P_{\dt}^{Z, Y, X, Y, Z} = \mathrm{e}^{\dt Z/2} \mathrm{e}^{\dt Y/2}  \mathrm{e}^{\dt X} \mathrm{e}^{\dt Y/2} \mathrm{e}^{\dt Z/2}.
\]
For instance, the second order scheme associated with~$P_{\dt}^{B, A, \gamma C, A, B}$ reads
\begin{equation}
  \label{eqn:BACAB}
  \left\{
    \begin{aligned}
      p^{n+\frac{1}{3}} &= p^n - \frac{\dt}{2} \nabla V(q^n) \COMMA \\
      q^{n+\frac{1}{2}} &= q^n + \frac{\dt}{2} M^{-1} p^{n + \frac{1}{3}} \COMMA \\
      p^{n+\frac{2}{3}} &= \alpha_{\dt} p^{n+\frac{1}{3}} + \sqrt{\frac{1-\alpha_{\dt}^2}{\beta} M} G^n\COMMA  \\
      q^{n+1} &= q^{n+\frac{1}{2}} + \frac{\dt}{2} M^{-1} p^{n+\frac{2}{3}} \COMMA \\
      p^{n+1} &= p^{n+\frac{2}{3}} - \frac{\dt}{2} \nabla V(q^{n+1}) \COMMA
    \end{aligned}
  \right.
 \end{equation}
where $(G^n)_{n \geq 0}$ are independent and identically distributed standard $D$-dimensional Gaussian random vectors;
while the second order scheme associated with~$P_{\dt}^{\gamma C,B, A, B,\gamma C}$ reads
\begin{equation}
  \label{eqn:CBABCB}
  \left\{
  \begin{aligned}
    p^{n+\frac{1}{4}} &= \alpha_{\dt/2} p^{n} + \sqrt{\frac{1-\alpha_{\dt/2}^2}{\beta} M} G_1^n\COMMA  \\
    p^{n+\frac{1}{2}} &= p^{n+\frac14} - \frac{\dt}{2} \nabla V(q^n) \COMMA \\
    q^{n+1} &= q^n + \dt M^{-1} p^{n + \frac{1}{2}} \COMMA \\
    p^{n+\frac{3}{4}} &= p^{n+\frac12} - \frac{\dt}{2} \nabla V(q^{n+1}) \COMMA \\
    p^{n+1} &= \alpha_{\dt/2} p^{n+\frac34} + \sqrt{\frac{1-\alpha_{\dt/2}^2}{\beta} M} G_2^n\COMMA  \\
    \end{aligned}
  \right.
 \end{equation}
where $(G_1^n)_{n \geq 0}$ and $(G_2^n)_{n \geq 0}$ are two independent families of independent and identically distributed standard $D$-dimensional Gaussian random vectors.

We recall in the next subsections various ergodicity results and error estimates for the schemes under consideration.

\subsubsection{Ergodicity results.}
One preliminary result used to prove the ergodicity of the numerical schemes is that the functions~$\mKS$ are Lyapunov functions for the numerical schemes under consideration (see~\cite[Lemma~2.7]{leimkuhler2016computation}): for any~$s^* \in \mathbb{N}$, there exists constants $\lambda>0$ and $K \in \mathbb{R}_+$ and a timestep~$\dt^*>0$ such that, for all $s = 1, 2, \ldots, s^*$ and $0 < \dt \leq \dt^*$,
\[
\forall n \in \mathbb{N}, \qquad (P_{\dt}^n \mKS)(y) \leq \mathrm{e}^{-\lambda n \dt} \mKS(y) + K \PERIOD
\]
In particular, there exists~$C \in \mathbb{R}_+$ such that
\begin{equation}
  \label{eqn:K-evolution-estimate}
  \forall n \in \mathbb{N}, \qquad (P_{\dt}^n \mKS)(y) \leq C \mKS(y) \PERIOD
\end{equation}
With this result at hand, the following theorem regarding the ergodicity of the associated Markov chain $\{y^n\}_{n\geq 0}$ is obtained from the results in~\cite{leimkuhler2016computation,EDMS21}.

\begin{theorem}
  \label{thm:ergodicity-splitting}
  Consider $s^* \geq 1$. For any $0 < \gamma < \infty$, there exists $\dt^* > 0$ such that, for any $0 < \dt \leq \dt^*$, the Markov chain associated with a first or second order splitting scheme has a unique invariant probability measure $\pi_{\gamma, \dt}$, which integrates the scale functions~$\mKS$ for all $s = 1, 2, \ldots, s^*$ uniformly in $\dt \in (0,\dt^*]$:
  \begin{equation}
    \label{eqn:integrate-scale-function}
    \forall s \in \{1, 2, \ldots, s^*\}, \qquad \sup_{0 < \dt \leq \dt^*}\int_{\mY} \mKS \, d\pi_{\gamma, \dt} < \infty.
  \end{equation}
  Moreover, there exist constants $\lambda = \lambda(s^*, \gamma) > 0$ and $K = K(s^*, \gamma) > 0$ such that, for any observable~$f \in B_{\mKS}^{\infty}(\mY)$,
  \begin{equation*}
    \forall n \in \mathbb{N},
    \qquad
    \left\| P_{\dt}^n f - \int_{\mY} f \,d\pi_{\gamma, \dt}\right\|_{B^\infty_{\mKS}} \leq K \mathrm{e}^{-\lambda n \dt} \|f\|_{B_{\mKS}^{\infty}} \PERIOD
  \end{equation*}
\end{theorem}

An important consequence of the above ergodicity result is the uniform control of the solution to the discrete time Poisson equation
\begin{equation}
  \label{eqn:disc-Poisson-eqn}
  - \LOPERH\WHF = f - \pi_{\gamma,\dt}(f)\COMMA
\end{equation}
where $\LOPERH$ is defined in~\eqref{eqn:discrete_generator}. Upon defining the subspace of~$B_{\mKS}^{\infty}(\mY)$ of functions with average~0 with respect to the probability measure~$\pi_{\gamma, \dt}$:
\[
B_{\mKS,\dt}^{\infty}(\mY) = \left\{f \in B_{\mKS}^{\infty}(\mY) ~\left|~ \int_{\mY} f \, d\pi_{\gamma, \dt} = 0\right. \right\}\COMMA
\]
we can recall the following result (see~\cite{EDMS21} and~\cite[Corollary~2.10]{leimkuhler2016computation}).

\begin{corollary}
  \label{cor:bounded-disc-Poission-solution}
  Consider an integer $s^* \geq 1$ and $\gamma\in (0,\infty)$. There exists $\dt^* > 0$ such that, for all $ s \in \{ 1, \ldots, s^*\}$,
  \begin{equation}
    \label{eqn:bounded-disc-resolvent}
    \sup_{0<\dt \leq \dt^*}\left\|\WHF\right\|_{B_{\mKS}^{\infty}} \leq \frac{2K}{\lambda}\left\|f \right\|_{B_{\mKS}^{\infty}} \COMMA
\end{equation}
with the same constants $K,\lambda$ as in Theorem~\ref{thm:ergodicity-splitting}.
\end{corollary}

\subsubsection{Error estimates on average properties.}
In the remainder of this work, we assume that the initial configuration $y^0 = (q^0, p^0)$ follows a distribution~$\mu_0$ which integrates all scale functions~$\mKS$, namely
\[
\forall s \geq 1, \qquad \int_{\mY}  \mKS \,d\mu_0 < \infty\PERIOD
\]
To simplify the notation, we denote by $\bE_{\dt} = \bE_{\dt}^{\mu_0}$ the expectation taken with respect to all initial conditions $y^0 \sim \mu_0$ and for all realizations of the Markov chain $\{y^n\}_{n\geq 0}$ starting from~$y^0$. 

The following result quantifies two sources of bias in the estimation of averages with respect to the target measure~$\pi$ by averages over one realization of the Markov chain: (i) a systematic bias arising from the use of a finite timestep~$\dt$; (ii) a truncation bias from the finite integration time. The former bias has been carefully studied in~\cite{leimkuhler2016computation} for Langevin dynamics, following the framework initiated by~\cite{TT90,talay2002stochastic}; while the latter bias is a direct consequence of Theorem~\ref{thm:ergodicity-splitting} (see also~\cite{mattingly2010convergence}).

\begin{theorem}[Bias of splitting schemes]
  \label{thm:bias-splitting}
  Fix~$\gamma \in (0,\infty)$ and an observable~$f \in \mS$, and consider either any first order ($\alpha = 1$) or second order ($\alpha =2$) splitting scheme. There exist $\dt^* > 0$ and $C_f > 0$ (which both depend on~$f$) such that, for any $0 < \dt \leq \dt^*$ and integer $N \geq 1$, 
  \begin{equation}
    \label{eqn:bias-splitting}
    \left|\frac{1}{N}\sum_{n=0}^{N-1} \bE_{\dt}\{f(y^n)\} - \int_{\mY} f \,d\pi  \right| \leq  C_f \left( \frac{1}{N\dt} + \dt^{\alpha}  \right)\PERIOD
  \end{equation}
\end{theorem}

\begin{proof}
  Denote by~$s \in \mathbb{N}$ an integer such that $f \in B^\infty_{\mKS}(\mY)$, and consider the timestep~$\dt^*$ as given by Theorem~\ref{thm:ergodicity-splitting}. It is a direct consequence of Theorem~$2.13$ ($\alpha = 1$) and Theorem~$2.16$ ($\alpha = 2$) in~\cite{leimkuhler2016computation} that there exists a constant $C_f \in \mathbb{R}_+$ such that, for all $\dt \in (0,\dt^*]$, 
  \begin{equation}
    \label{eq:error_ergodic_avg_splitting}
\left|\int_{\mY} f \, d\pi_{\gamma, \dt} - \int_{\mY}  f \, d\pi \right| \leq C_f \dt^{\alpha} \PERIOD
\end{equation}
On the other hand, by Theorem~\ref{thm:ergodicity-splitting},
\begin{equation*}
  \begin{aligned}
  & \left|\frac{1}{N}\sum_{n=0}^{N-1} \bE_{\dt}\{f(y^n)\} - \int_{\mY} f \,d\pi_{\gamma, \dt}  \right|
  \leq \frac{1}{N}\sum_{n=0}^{N-1} \int_\mY \left|P_\dt^n f - \int_{\mY} f \,d\pi_{\gamma, \dt}\right| d\mu_0 \\
  & \qquad \leq \frac{K}{N} \mu_0(\mKS)\|f\|_{B_{\mKS}^{\infty}}\sum_{n=0}^{N-1} \mathrm{e}^{-\lambda n \dt}
  = \frac{K(1-\mathrm{e}^{-\lambda N\dt})}{N(1-\mathrm{e}^{\lambda \dt})} \mu_0(\mKS)\|f\|_{B_{\mKS}^{\infty}} 
  \leq \frac{2K}{\lambda N\dt} \mu_0(\mKS) \|f\|_{B_{\mKS}^{\infty}} \COMMA
  \end{aligned}
\end{equation*}
upon possibly reducing~$\dt^*$. The error estimate~\eqref{eqn:bias-splitting} then immediately follows from the triangle inequality, with $C = \max\{C_f, 2K \mu_0(\mKS)\|f\|_{B_{\mKS}^{\infty}}/\lambda\}$.
\end{proof}

\subsection{Approximate solutions to Poisson equations}
\label{subsec:approx-Poisson-equation}

The formal derivation of MP estimators in Section~\ref{sec:general-continuous-setting} crucially relies on the fact that solutions to the discrete Poisson equation~\eqref{eqn:disc-Poisson-eqn}, which belong to functional spaces such as~$B^\infty_{\mKS}(\mY)$ when~$f \in \mathcal{S}$, can in fact be approximated by smooth functions, as given by~\eqref{eqn:approx-Poisson-solution} below. As in~\cite[Section~4.3]{leimkuhler2016computation} (see also~\cite{lelievre2016partial,plechac2019convergence}), we rely on approximate inverse operators to this end. We briefly recall in this section how these operators are constructed, and provide some technical estimates on expectations involving approximate solutions to Poisson equations.

For the splitting schemes considered in Section~\ref{sec:splitting_Lang}, the evolution operator~$P_\dt$ can be expanded in powers of~$\dt$ as follows, for a given function $f\in \mathcal{S}$:
\begin{equation}
  \label{eq:expansion_P_dt}
  P_{\dt}f = f + \dt\AOP_1f + \ldots + \dt^{\alpha+1}\AOP_{\alpha+1}f + \dt^{\alpha+2} r_{\alpha,\dt,f}\COMMA
\end{equation}
with $\alpha = 1$ for first order and $\alpha = 2$ for second order schemes. The operators~$\AOP_i$ can be systematically identified by the Baker--Campbell--Hausdorff formula~\cite[Section~III.4.2]{hairer2006geometric}, while the functions~$r_{\alpha, \dt,f}$ are in~$\mathcal{S}$ and there exists~$s \in \mathbb{N}$, $\dt^*>0$ and~$K\in \mathbb{R}_+$ such that~$\|r_{\alpha, \dt,f}\|_{B^\infty_{\mKS}} \leq K<\infty$ for~$\dt \in (0,\dt^*]$. 
We next note that the discrete time Poisson equation~\eqref{eqn:disc-Poisson-eqn} can be rewritten as
\begin{equation}
  \label{eqn:discrete-Poisson-equation-projected}
  -\LOPERHPI \WHF = f - \pi(f)\COMMA
\end{equation}
with (recalling the definition~\eqref{eqn:projection-operator} for~$\Pi$)
\[
\LOPERHPI = \Pi\LOPERH\Pi \equiv \Pi\left[\frac{1}{\dt}(P_h - I)\right]\Pi\PERIOD
\]
In order to find an approximation to $\WHF$, we introduce~$\TBOP = \TAOP_2 + \dt \TAOP_3 + \ldots + \dt^{\alpha-1} \TAOP_{\alpha+1}$, where $\TAOP_j = \Pi \AOP_j \Pi$ for $j = 1, 2, \ldots, \alpha+1$ are operators mapping~$\mS_0$ to itself, so that $\LOPERHPI f =\TAOP_1 f + \dt \TBOP f + \dt^{\alpha+2} r_{\alpha,\dt,f}$. We next truncate the formal series expansion of the inverse 
\[
\left(\TAOP_1 + \dt \TBOP\right)^{-1} = 
\TAOP_1^{-1} - \dt \TAOP_1^{-1}\TBOP \TAOP_1^{-1} + \dt^2 \TAOP_1^{-1}\TBOP\TAOP_1^{-1}\TBOP\TAOP_1^{-1} + \ldots\PERIOD
\]
up to terms involving at most~$\alpha$ instances of~$\TBOP$. This motivates introducing
\[
\TQOP = \TAOP_1^{-1}\sum_{j=0}^{\alpha} (-1)^j \dt^j \left(\TBOP\TAOP_1^{-1}\right)^j\COMMA
\]
which satisfies the following identity on~$\mathcal{S}_0$:
\[
\LOPERHPI \TQOP f = \left(\TAOP_1 + \dt \TBOP\right)\TQOP f + \dt^{\alpha+2} r_{\alpha,\dt,\TQOP f} = \Pi f + (-1)^{\alpha} \dt^{\alpha+1} \left(\TBOP\TAOP_1^{-1}\right)^{\alpha+1} f + \dt^{\alpha+2} \widetilde{r}_{\alpha,\dt,f} \PERIOD
\]
Note that the functions~$\widetilde{r}_{\alpha, \dt, f}$ are in~$\mathcal{S}$ and there exists~$s \in \mathbb{N}$, $\dt^*>0$ and~$K \in \mathbb{R}_+$ such that~$\|\widetilde{r}_{\alpha, \dt, f}\|_{B^\infty_{\mKS}} \leq K<\infty$ for any~$\dt \in (0,\dt^*]$. This can be seen from the fact that $\widetilde{r}_{\alpha, \dt,f} = r_{\alpha, \dt, \TQOP f}$ is a finite sum of terms of the form~$\dt^j r_{\alpha, \dt, g_j}$ with~$g_j \in \mS$.
We finally define the {\it approximate inverse operator} $Q_{\dt}$ by substituting the formula for~$\TBOP$ into $\TQOP $ and keeping only the terms of order at most~$\dt^{\alpha}$, \emph{i.e.} 
\[
Q_{\dt} \triangleq \mathcal{Q}_0 + \dt \mathcal{Q}_1 + \ldots + \dt^{\alpha-1}\mathcal{Q}_{\alpha-1} + \dt^{\alpha} \mathcal{Q}_{\alpha}\COMMA
\] 
where the operators~$\mathcal{Q}_{j}$ for $0 \leq j \leq \alpha$ are operators mapping~$\mathcal{S}_0$ to itself. Therefore, $Q_{\dt}$ also maps~$\mS_0$ to itself.

\begin{remark}
\label{rmk:second_order_Q1}
To derive MP estimators in Section~\ref{sec:MP_for_Langevin} (in particular, the second order MP estimator in Section~\ref{subsec:2nd-order-MP}), we only need the explicit expressions of $\mathcal{Q}_0 = \TAOP_1^{-1}$ and $\mathcal{Q}_1 = -\TAOP_1^{-1} \TAOP_2 \TAOP_1^{-1}$. For schemes of weak order~2, it holds~$\AOP_1 = \mL_{\gamma}$ and $\AOP_2 = \mL_{\gamma}^2/2$, so that~$\mathcal{Q}_1 = -\Pi/2$.
\end{remark}

We are now in a position to define the {\it approximate solution to the discrete time Poisson equation} 
\begin{equation}
  \label{eqn:approx-Poisson-solution}
  \WTF = -Q_{\dt}(f - \pi(f))\PERIOD 
\end{equation}
We emphasize that the function~$\WTF$ indeed belongs to~$\mathcal{S}_0$ when~$f \in \mS$, while the solution~$\WHF$ to the Poisson equation~\eqref{eqn:disc-Poisson-eqn-general} does not. Applying the operator $-\LOPERHPI$ 
to both sides of~\eqref{eqn:approx-Poisson-solution} leads to 
\begin{equation}
  \label{eqn:phi}
  -\LOPERHPI \WTF = f - \pi(f) + \dt^{\alpha+1} \phi_{\alpha, \dt, f}
\end{equation}
for some function $\phi_{\alpha, \dt, f} \in \mathcal{S}_0 $, and there exists~$s \in \mathbb{N}$, $\dt^*>0$ and~$K\in \mathbb{R}_+$ such that
\begin{equation}
  \label{eqn:control-of-phi}
  \forall \dt \in (0,\dt^*], \qquad \|\phi_{\alpha, \dt,f}\|_{B^\infty_{\mKS}} \leq K \PERIOD
\end{equation}
Moreover, in view of~\eqref{eqn:K-evolution-estimate}, for any $a^* \in \mathbb{N}$, upon possibly increasing~$s$ and~$K$ and decreasing~$\dt^*$,
\begin{equation}
 \label{eq:second-moment-bound-widetilde}
 \forall a \in \{1, \ldots, a^*\},\qquad  \forall \dt \in (0,\dt^*], \qquad \forall n \geq 0, \qquad \bE_{\dt}\left\{ \left|\WTF(y^n)\right|^{2a} \right\} \leq K \left\|f\right\|_{B_{\mKS}^{\infty}}^{2a} \PERIOD
\end{equation}

Finally, from the definitions~\eqref{eqn:cont-Poisson-eqn} and~\eqref{eqn:approx-Poisson-solution}, we can compare~$\WTF$ with~$\widehat{f}$ when $\AOP_1 = \mL_{\gamma}$ (\emph{i.e.} the numerical scheme is at least of weak order~1):
\begin{equation}
  \label{eqn:diff-discrete-continuous-Poisson-pointwise}
  \WTF-\widehat{f} = -\left[\dt  \mathcal{Q}_1 + \dt^2 \mathcal{Q}_2 + \ldots + \dt^{\alpha-1}\mathcal{Q}_{p-1} + \dt^{\alpha} \mathcal{Q}_{\alpha}\right](f - \pi(f))\PERIOD
\end{equation}
In view of Remark~\ref{rmk:second_order_Q1}, this corresponds to the equalities~\eqref{eqn:approx-operator-preview} in the general argument of Section~\ref{sec:general-continuous-setting}.

\section{Martingale product estimators for Langevin dynamics}
\label{sec:MP_for_Langevin}

We are now able to derive the MP estimators for the linear response of Langevin dynamics, by following the general strategy presented in Section~\ref{sec:general-continuous-setting}. We first discuss in Section~\ref{sec:choice_g_Lang} the choice of the function~$g$ appearing in the MP estimator, which allows us to study the consistency of both the first and second order MP estimators in Sections~\ref{subsec:1st-order-MP} and~\ref{subsec:2nd-order-MP}, respectively. The boundedness of the variance of MP estimators is proved in Section~\ref{subsec:MP-var-analysis}. In all this section, we rewrite the elementary evolutions associated with first or second order splittings as
\[
  y^{n+1} = y^n + \Phi_{\dt}(y^n; G^n),
  \qquad
  \Phi_{\dt}(y; G) = \begin{pmatrix} \Phi_{\dt}^q(y; G) \\ \Phi_{\dt}^p(y; G) \end{pmatrix} \COMMA
\]
where the expression of the increment function~$\Phi_{\dt}$ depends on the specific scheme under consideration, and $(G^n)_{n \geq 0}$ are independent and identically distributed standard Gaussian random vectors. It may be the case that two independent Gaussian vectors $G^{n,1},G^{n,2}$ are needed per timestep for second order splittings such as~$\mathrm{e}^{\dt C/2} \mathrm{e}^{\dt A/2} \mathrm{e}^{\dt B} \mathrm{e}^{\dt A/2} \mathrm{e}^{\dt C/2}$, in which case~$G^n=(G^{n,1},G^{n,2})$.

\subsection{Choice of the function~$g$ in the discrete auxiliary martingale}
\label{sec:choice_g_Lang}

Given our specific choice for the external forcing $F(q)$, the
equation~\eqref{eqn:constraint-g} on the function $g$ which appears 
in the MP estimator is easily solved by choosing
\begin{equation}
    \label{eq:solutio_g_MP_Langevin}
    g(q,p) = \frac{\beta}{2\gamma} p^\T F(q).
\end{equation}

The solution is not unique\footnote{We do not make use of this flexibility here, but an interesting question is whether one of the possible solutions leads to a final estimator with minimal variance.} since any function of the~$q$ variable only can be added to~$g$. The important point here is that we have an explicit solution which can be used to construct the numerical scheme. 

\begin{remark}\label{rem:gfun}
For more general forcings depending both on momenta and positions, i.e.,
forces $F(q,p)$, a solution to~\eqref{eqn:constraint-g} is given by 
\begin{equation}\label{eqn:gfun}
g(q,p) = \frac{\beta}{2\gamma} \left(\nabla_p^* \nabla_p \right)^{-1} \nabla_p^* F(q,p)\,,
\end{equation}
where $\nabla_p^*$ is defined in~\eqref{eq:PDE_g}.
This function is well-defined $q$ by $q$, since $\nabla_p^* F$ has zero average with respect to the marginal measure~$\kappa(dp)$ of~$\pi$ in the $p$-variable (\emph{i.e.}, $\kappa$ is the Gaussian measure in~$p$ with mean zero and covariance~$M\beta$), and $\nabla_p^* \nabla_p$ is coercive on the subspace of functions in~$L^2(\kappa)$ with mean zero with respect to~$\kappa$ (thanks to the Poincar\'e inequality satisfied by~$\kappa$). On the other hand, the expression for~$g$ is not explicit in general. This prevents explicitly constructing the discrete martingale by expansions in powers of~$\dt$ of~$\DMN(g)$. When the expression of~$g$ is explicit (which is the case when $F(q,p)$ depends polynomially in~$p$, for example), the methodology presented in this section can be adapted in a straightforward way, see Appendix~\ref{app:secondorder} for more details.
\end{remark}

\subsection{First order martingale product estimator}
\label{subsec:1st-order-MP}

A simple computation shows that, for all the first order splitting schemes,
\begin{equation}
\label{eq:expansion_first_order_splitting}
  \Phi_{\dt}^p(y^n; G^n) = \dt^{1/2}\sqrt{\frac{2\gamma}{\beta}} G^n + \Delta_{\dt}(y^n; G^n),
\end{equation}
where the remainder term~$\Delta_{\dt}(y^n; G^n)$ and~$\Phi_\dt^q$ are of order~$\dt$, in the following sense: for any~$a^* \in \mathbb{N}$, there exists~$C_* \in \mathbb{R}_+$ and~$\dt^* >0$ such that 
\begin{equation}
  \label{eq:remainders_first_order_splitting}
  \forall a \in \{1,\ldots,a^*\}, \quad \forall \dt \in (0,\dt^*],
  \qquad
  \bE_{\dt}\{ |\Delta_{\dt}(y^n;G^n)|^a\} + \bE_{\dt}\{ |\Phi_\dt^q(y^n;G^n)|^a\} \leq C_* \dt^{a}.
\end{equation}
As heuristically derived in Section~\ref{sec:first_order_general}, the first order MP estimator based on an arbitrary first order splitting scheme is
\begin{equation}
  \label{eqn:first-order-MP-estimator}
  \mathcal{E}_{\dt, N}^{\mathrm{MP}1}(f) =\frac{1}{N} \sum_{n=0}^{N-1}\left(f(y^n) - \pi_{\gamma, \dt}(f)\right) z^N\COMMA
\end{equation}
where 
\begin{equation}
  \label{eqn:1st-order-weight}
  z^N = \dt^{1/2} \sqrt{\frac{\beta}{2\gamma}} \sum_{n=0}^{N-1} F(q^n)^{\T} G^n
\end{equation}
is the auxiliary discrete martingale, obtained from expanding $g(y^{n+1}) - (P_h g)(y^n)$ up to order $\sqrt{h}$ in~\eqref{eqn:general-first-order-MP} and replacing~$g$ by its expression~\eqref{eq:solutio_g_MP_Langevin}. Note that the weight process remains the same for all the first order splittings considered in Section~\ref{sec:splitting_Lang}. The first order MP estimator coincides with the CLR estimator (see~\cite{plechac2019convergence} for the overdamped Langevin setting). The benefit of the proposed martingale product derivation is that it provides a systematic approach for discovering higher order MP schemes, as discussed in Section~\ref{subsec:2nd-order-MP}. The reason why we carefully study first order estimators is that the proof of their consistency paves the way for the proof of consistency of second order estimators.

The main result concerning the consistency of first order MP estimators is the following.

\begin{theorem}[Bias of first order MP estimators]
  \label{thm:1st-order-bias}
  Consider any of the first order splitting schemes introduced in Section~\ref{sec:splitting_Lang}, and fix an observable~$f \in \mS$. Then there exist $\dt^* > 0$ and $C \in \mathbb{R}_+$ such that, for any $0<\dt \leq \dt^*$ and $N \geq 1$ satisfying~$N\dt \geq 1$,
\begin{equation}\label{eqn:1st-order-MP-bias}
\left|\bE_{\dt}\left\{ \mathcal{E}_{\dt, N}^{\mathrm{MP}1}(f)\right\} - \LIN (f)\right| \leq C\left(\dt + \frac{1}{\sqrt{N\dt}}\right)\PERIOD
\end{equation}
\end{theorem}

\begin{proof}
  Throughout the proof, we denote by~$C$ a positive constant that may vary line by line (but does not depend on~$N$ and~$\dt$) and by~$s$ a sufficiently large positive integer. We also repeatedly use the following bound, directly obtained from~\eqref{eqn:K-evolution-estimate}: for given functions~$\varphi_1,\varphi_2 \in \mathcal{S}$ (so that $\varphi_1 \varphi_2 \in \mS$), there exist~$K \in \mathbb{R}_+$ and~$\dt > 0$ such that
  \begin{equation}
    \label{eq:bound_sum_expectation_functions}
    \forall \dt \in (0,\dt^*], \quad \forall N \geq 1, \qquad \frac{1}{N} \left| \sum_{n=0}^{N-1} \bE_{\dt}\left\{ \varphi_1(y^n) \varphi_2(y^n) \right\}\right| \leq K.
  \end{equation}

  We start by rewriting the MP estimator in a form more amenable for our estimates. Since $\bE_{\dt}\{z^N\} = 0$, one can easily verify that
\[
\bE_{\dt}\left\{ \left(f(y^n) - \pi_{\gamma, \dt}(f)\right)  z^N  \right\} = \bE_{\dt}\left\{  \left(f(y^n) - \pi(f)\right) z^N \right\}.
\]
In view of~\eqref{eqn:phi} (with~$\alpha=1$), it holds
\begin{equation}
  \label{eqn:1st-order-estimate}
  \begin{aligned}
    & \bE_{\dt}\left\{ \mathcal{E}_{\dt, N}^{\mathrm{MP}1}(f)\right\} = \frac{1}{N\dt} \sum_{n=0}^{N-1} \bE_{\dt}\left\{  \left(\Pi (I - P_{\dt}) \Pi \WTF(y^n)\right) z^N  \right\} - \frac{\dt^2}{N} \sum_{n=0}^{N-1} \bE_{\dt}\left\{ \phi_{1, \dt, f}(y^n) z^N \right\} \\
    & \qquad = \frac{1}{N\dt} \sum_{n=0}^{N-1} \bE_{\dt}\left\{ \left( \WTF(y^{n+1}) - (P_{\dt} \WTF)(y^n)   \right)z^N \right\} + \frac{1}{N\dt} \bE_{\dt}\left\{ \left( \WTF(y^0) - \WTF(y^N)  \right) z^N \right\}\\
    & \qquad \ \ \ - \frac{\dt^2}{N} \sum_{n=0}^{N-1} \bE_{\dt}\left\{ \phi_{1, \dt, f}(y^n) z^N \right\}\PERIOD
  \end{aligned}
\end{equation}
Let us first show that the last two terms in the last right hand side vanish as~$\dt \to 0$. Since the functions~$\phi_{1, \dt, f}$ belong to some space~$B_{\mKS}^{\infty}(\mY)$ for~$\dt$ sufficiently small by~\eqref{eqn:control-of-phi}, the estimates in Proposition~\ref{prop:elementary-term} immediately imply that
\begin{equation}
  \label{eqn:1st-order-estimate-1}
  \frac{\dt^2}{N} \left| \sum_{n=0}^{N-1} \bE_{\dt}\left\{ \phi_{1, \dt, f}(y^n) z^N \right\} \right| \leq C\dt^{3/2}
\end{equation}
for $\dt$ sufficiently small. For the second term in the last right hand side of~\eqref{eqn:1st-order-estimate}, we use~\eqref{eq:bound_sum_expectation_functions} and the assumption that each component of $F$ is a smooth function to write
\begin{equation}
  \label{eq:bound_zN^2}
  \bE_{\dt}\left\{\left(z^N\right)^2\right\} = \frac{\beta\dt}{2\gamma}\sum_{n=0}^{N-1} \bE_\dt\left\{ \left|F(q^n)\right|^2 \right\}\leq C N\dt \PERIOD
\end{equation}
A Cauchy--Schwarz inequality together with~\eqref{eq:second-moment-bound-widetilde} then leads to
\begin{equation}\label{eqn:1st-order-estimate-2}
\frac{1}{N\dt} \left| \bE_{\dt}\left\{ \left( \WTF(y^0) - \WTF(y^N)  \right) z^N \right\} \right| \leq \frac{C}{\sqrt{N\dt}} \left\|f \right\|_{B_{\mKS}^{\infty}},
\end{equation}
which contributes to the finite integration time error~$C /\sqrt{N\dt}$ in the bias.

It only remains to estimate the difference between the first term on the right hand side of~\eqref{eqn:1st-order-estimate} and~$\LIN(f)$. Since the martingale increments~$z^{m+1}-z^m$ and~$\DMN(\WTF)$ are conditionally independent for $m \neq n$ (recalling the notation in~\eqref{eqn:obs_increment}), 
\[
\sum_{n=0}^{N-1} \bE_{\dt}\left\{ \left( \WTF(y^{n+1}) - (P_{\dt} \WTF)(y^n)   \right) z^N \right\}
 =
\sum_{n=0}^{N-1} \bE_{\dt}\left\{ \DMN(\WTF)(z^{n+1} - z^n) \right\}.
\]
Now, by a Taylor expansion and in view of~\eqref{eq:remainders_first_order_splitting}, 
\begin{equation}
  \label{eq:Taylor_WTF}
  \DMN(\WTF) = \sqrt{\frac{2\gamma\dt}{\beta}} \left(G^n\right)^\T \nabla_p \WTF(y^n) + \dt \psi(y^n;G^n) + \dt^{3/2}\Upsilon_\dt(y^n;G^n) \COMMA
\end{equation}
with $\bE_\dt(\psi (y^n;G^n) G^n) = 0$ and terms~$\psi(y^n;G^n),\Upsilon_\dt(y^n;G^n)$ of order~1 in the following sense: for any~$a \geq 1$, there exist~$C \in \mathbb{R}_+$ and~$\dt^*>0$ such that, for any~$0 \leq n \leq N$ and~$\dt \in (0,\dt^*]$, 
\begin{equation}
\label{eq:moment_bounds_psi_Upsilon}
    \bE_{\dt}\left\{ |\psi (y^n;G^n)|^a \right\} \leq C, 
    \qquad 
    \bE_{\dt}\left\{ |\Upsilon_\dt (y^n;G^n)|^a \right\} \leq C.
\end{equation}
    A Cauchy--Schwarz inequality then gives
\[
  \left|\sum_{n=0}^{N-1} \bE_{\dt}\left\{ \Upsilon_\dt (y^n;G^n)(z^{n+1} - z^n) \right\}\right| \leq \left(\sum_{n=0}^{N-1} \bE_{\dt}\left\{ \Upsilon_\dt (y^n;G^n)^2 \right\}\right)^{1/2}\left( \sum_{n=0}^{N-1} \bE_{\dt}\left\{ (z^{n+1}-z^n)^2\right\}\right)^{1/2},
\]
so that, since the second sum on the right hand side of the above inequality is bounded by~$CN\dt$,
\begin{equation}
  \label{eq:1st-order-bound_remainder_martingale}
  \frac{1}{N\dt} \left|\sum_{n=0}^{N-1} \bE_{\dt}\left\{ \left[ \dt \psi(y^n;G^n)+\dt^{3/2} \Upsilon_\dt (y^n;G^n)\right](z^{n+1} - z^n) \right\}\right| \leq C \dt.
\end{equation}
Moreover, by~\eqref{eqn:diff-discrete-continuous-Poisson-pointwise} with~$\alpha=1$,
\[
  \begin{aligned}
    & \frac{1}{N\dt} \sqrt{\frac{2\gamma\dt}{\beta}} \sum_{n=0}^{N-1} \bE_{\dt}\left\{ \left(G^n\right)^\T \nabla_p \WTF(y^n) (z^{n+1} - z^n) \right\} = \frac{1}{N}\sum_{n=0}^{N-1} \bE_{\dt}\left\{ F(y^n)^\T \nabla_p \WTF(y^n) \right\} \\
    & \qquad = \frac{1}{N}\sum_{n=0}^{N-1} \bE_{\dt}\left\{ F(y^n)^\T \nabla_p \widehat{f}(y^n) \right\} - \frac{\dt}{N}\sum_{n=0}^{N-1} \bE_{\dt}\left\{ F(y^n)^\T \nabla_p \mathcal{Q}_1 \Pi f(y^n) \right\} \PERIOD    
  \end{aligned}
\]
The second sum in the previous equality is uniformly bounded by~$C\dt$ in view of~\eqref{eq:bound_sum_expectation_functions}, while we use Theorem~\ref{thm:bias-splitting} for the first one: 
\begin{equation}
  \label{eq:1st-order_estimates_dominant_term}
  \left| \frac{1}{N}\sum_{n=0}^{N-1} \bE_{\dt}\left\{ F(y^n)^\T \nabla_p \widehat{f}(y^n) \right\} - \int_{\mY} F^\T \nabla_p\widehat{f} \,d\pi \right| \leq C \left( \dt + \frac{1}{N\dt} \right). 
\end{equation}
The integral on the left-hand side of the previous equality is equal to~$\LIN(f)$ by the reformulation~\eqref{eqn:Langevin-lin-res-reformulation} for the linear response. The desired bias estimate finally follows by combining~\eqref{eqn:1st-order-estimate} and the estimates provided by~\eqref{eqn:1st-order-estimate-1}, \eqref{eqn:1st-order-estimate-2}, \eqref{eq:1st-order-bound_remainder_martingale} and~\eqref{eq:1st-order_estimates_dominant_term}.
\end{proof}

\subsection{Second order martingale product estimator}
\label{subsec:2nd-order-MP}

We present in this section second order MP estimators based on second order splittings. For first order MP estimators, the auxiliary discrete martingales~$\{z^n\}_{n \geq 0}$ are identical for different first order splittings. However, this is not the case for second order MP estimators. We choose to first present in Section~\ref{sec:general_2nd_order} a general result on second order MP estimators, following the strategy outlined in Section~\ref{sec:construction_discrete_martingale}. The so-obtained estimators are however not explicit, which is why we make them precise in Section~\ref{sec:specific_2nd_order} for the various splitting schemes we consider. 

\subsubsection{General second order estimators.}
\label{sec:general_2nd_order}

The second order MP estimator obtained from~\eqref{eqn:general-second-order-MP} is
\begin{equation}
  \mathcal{E}_{\dt, N}^{\mathrm{MP}2} = \frac{1}{N}\sum_{n=0}^{N-1}\left[ f(y^n) -\pi_{\gamma,\dt}(f)\right] \DMN\LGC\,,
   \label{eq:MP2_estimator}
\end{equation}
where we recall
\[
\DMN\LGC = \LGC(y^{n+1}) - P_{\dt} \LGC (y^n)\PERIOD\notag
\]
The second order consistency of this estimator is rigorously established by the following result.

\begin{theorem}[Bias of second order MP estimators]
  \label{thm:2nd-order-bias}
  Consider any of the second order splitting schemes introduced in Section~\ref{sec:LD-application}, and fix an observable $f \in \mS$. Then there exist $\dt^* > 0$ and $C \in \mathbb{R}_+$ such that, for any $0<\dt \leq \dt^*$ and $N \geq 1$,
  \begin{equation}
    \label{eqn:2nd-order-MP-bias}
    \left|\bE_{\dt}\left\{ \mathcal{E}_{\dt, N}^{\mathrm{MP}2}(f)\right\} - \LIN (f)\right| \leq C\left(\dt^2 + \frac{1}{\sqrt{N\dt}}\right)\PERIOD
  \end{equation}
\end{theorem}

\begin{proof}
  Throughout the proof, we denote by~$C$ a positive constant and by~$s$ a positive integer that is sufficiently large. These two numbers may vary line by line. We introduce the auxiliary discrete martingale
  \[
    z^N = \sum_{n=0}^{N-1} \DMN\LGC.
  \]
  First note that~\eqref{eq:expansion_first_order_splitting} and~\eqref{eq:remainders_first_order_splitting} hold true for all second order splitting schemes as well. 
  Furthermore, from~\eqref{eqn:K-evolution-estimate} and  
  together with the fact that~$g,\mL g \in \mS$,
  \[
    \bE_\dt\left\{ \left(z^{n+1}-z^n\right)^2 \right\}  
    \leq 2\bE_\dt\left\{ \DMN(g)^2 \right\} + \frac{\dt^2}{2}\bE_\dt\left\{ \DMN(\mL g)^2 \right\}  \leq D \dt\COMMA
  \]
  for some constant $D> 0$. 
  Similarly to~\eqref{eqn:1st-order-estimate}, using~\eqref{eqn:phi} with~$\alpha=2$ and the conditional independence of the martingale increments $z^{m+1}-z^m$ and~$\WTF(y^{n+1}) - (P_{\dt} \WTF)(y^n)$ when $m \neq n$,
\begin{equation*}
\begin{split}
\bE_{\dt}\left\{ \mathcal{E}_{\dt, N}^{\mathrm{MP}2}(f)\right\} & = \frac{1}{N\dt} \sum_{n=0}^{N-1} \bE_{\dt}\left\{ \DMN(\WTF)(z^{n+1}-z^n) \right\}\\
& \ \ \ \ +
\frac{1}{N\dt} \bE_{\dt}\left\{ \left( \WTF(y^0) - \WTF(y^n)  \right) z^N \right\} - \frac{\dt^3}{N} \sum_{n=0}^{N-1} \bE_{\dt}\left\{ \phi_{2, \dt, f}(y^n) z^N \right\} \PERIOD
\end{split}
\end{equation*}
The second and third terms on the right hand side of the above equality can be uniformly controlled as in the proof of Theorem~\ref{thm:1st-order-bias}. Now, in view of~\eqref{eq:crucial_identity}, 
\[
  \bE_{\dt}\left\{ \DMN(\WTF)(z^{n+1}-z^n) \right\} %
  = \dt \bE_\dt\left\{ \Gamma_\dt\left(\WTF,g-\frac{\dt}{2}\mL g\right)(y^n)\right\},
\]
so that, with~\eqref{eqn:bias-splitting},
\[
\left| \frac{1}{N\dt} \sum_{n=0}^{N-1} \bE_{\dt}\left\{ \DMN(\WTF)(z^{n+1}-z^n) \right\} - \int_\mY \Gamma_\dt\left(\WTF,g-\frac{\dt}{2}\mL g\right) d\pi \right| \leq C \left(\dt^2 + \frac{1}{N\dt}\right) \COMMA
\]
where the constant~$C$ depends on~$F$ through~$g$, and on the observable~$f$ under consideration. At this stage, it only remains to show that the difference between~$\LIN(f)$ and the integral on the left hand side of the previous inequality is of order~$\dt^2$. To this end, we first use~\eqref{eq:expansion_P_dt} (which allows us to give a rigorous meaning to~\eqref{eqn:order2}) and~\eqref{eqn:diff-discrete-continuous-Poisson-pointwise} (with~$\mathcal{Q}_1\varphi = -\Pi\varphi/2$ for second order splitting schemes) to write
\[
  \Gamma_\dt\left(\WTF,g-\frac{\dt}{2}\mL g\right) = \Gamma\left(\widehat{f},g\right) + \dt\left[\Lambda\left(\widehat{f},g\right)- \frac12 \Gamma(\mL \widehat{f},g) - \frac12 \Gamma\left(\widehat{f},\mL g\right) \right] + \dt^2 r_{\dt,f} \COMMA
\]
with $\|r_{\dt,f}\|_{B^\infty_{\mKS}} \leq K$ for~$\dt$ sufficiently small. The result then follows from~\eqref{eqn:integrate-scale-function} and the fact that $\Gamma\left(\widehat{f},g\right) = F^\T \nabla_p \widehat{f}$, while the order~$\dt$ term in the expansion of~$\Gamma_\dt\left(\WTF,g-\dt\mL g/2\right)$ in powers of~$h$ has average~0 with respect to~$\pi$ in view of~\eqref{eq:int_Lambda_pi}.
\end{proof}

\begin{remark}
A closer inspection of the proof of Theorem~\ref{thm:2nd-order-bias} shows that similar results can be obtained for other discretization schemes as long as $\mA_1 = \mL_{\gamma}$ and~$\mA_2$ is proportional, but not necessarily equal, to~$\mL_{\gamma}^2$; see~\cite{FS17} for examples of such schemes, based on a Barker acceptance/rejection rule together with a second order discretization of overdamped Langevin dynamics. 
\end{remark}

\subsubsection{Application to second order splitting schemes.}
\label{sec:specific_2nd_order}
We now provide a more explicit expression for the estimator~\eqref{eq:MP2_estimator} depending on the numerical scheme at hand. We illustrate our approach for the second order splitting~$P_{\dt}^{B,A,\gamma C, A, B}$ considered in~\eqref{eqn:BACAB} and then discuss the extension to other second order splittings. The strategy is to expand the martingale increments~$\DMN(g-\dt\mL g/2)$ up to order~$\dt^{3/2}$ as
\begin{equation}
  \label{eq:expansion_martingale_increment_2nd_order}
  \begin{aligned}
    \DMN\left(g-\frac{\dt}{2}\mL g\right) = h^{1/2}\sqrt{\frac{\beta}{2\gamma}} F(y^n)^{\T} G^n & + \dt \psi_1(y^n,G^n) + \dt^{3/2}\psi_{3/2}(y^n,G^n) \\
    & + \dt^2 \psi_2(y^n,G^n) + \dt^{5/2} \Upsilon_\dt(y^n,G^n),
  \end{aligned}
\end{equation}
where the dominant term on the right hand side of course agrees with the discrete martingale increments in~\eqref{eqn:1st-order-weight} for the first order MP estimator. It then suffices to replace the last factor in~\eqref{eq:MP2_estimator} by the first three terms above. The remainder terms involving~$\psi_2$ and~$\Upsilon_\dt$ which arise when quantifying the bias of the estimator can be uniformly controlled by noting that $\bE_\dt(\psi_2(y^n,G^n) G^n) = 0$ and using uniform estimates on moments of~$\psi_2$ and~$\Upsilon_\dt$. We do not explicitly perform these manipulations as they follow very closely estimates used in the proofs of Theorems~\ref{thm:1st-order-bias} and~\ref{thm:2nd-order-bias}. 

The first step towards an equality such as~\eqref{eq:expansion_martingale_increment_2nd_order} is to expand the increment function of the scheme under consideration up to the order of $\dt^{3/2}$. For~\eqref{eqn:BACAB},   
\begin{equation}
\label{eq:expansion_BACAB}
\begin{split}
\Phi_{\dt}^q(y^n; G^n) &= \dt M^{-1} p^n + \dt^{3/2} \sqrt{\frac{\gamma}{2\beta}}M^{-1} G^n + \Delta_{\dt}^q(y^n; G^n),\\
\Phi_{\dt}^p(y^n; G^n) 
&= 
h^{1/2}\sqrt{\frac{2\gamma}{\beta}} G^n - \dt (\nabla V(q^n) + \gamma M^{-1} p^n) - \dt^{3/2}\sqrt{\frac{\gamma}{2\beta}}\gamma M^{-1} G^n
+ \Delta_{\dt}^p(y^n; G^n),
\end{split}
\end{equation}
with remainders terms of order~$\dt^2$: for any $a^* \in \mathbb{N}$, there exists~$C_* \in \mathbb{R}_+$ and~$\dt^* >0$ such that 
\[
  \forall a \in \{1,\ldots,a^*\}, \quad \forall \dt \in (0,\dt^*],
  \qquad
  \bE_{\dt}\{ |\Delta_{\dt}^q(y^n;G^n)|^a\} + \bE_{\dt}\{ |\Delta_{\dt}^p(y^n;G^n)|^a\} \leq C_* \dt^{2a} \PERIOD
\]
This leads to the following second order MP estimator for the scheme~\eqref{eqn:BACAB}:
\begin{equation}
  \label{eqn:second-order-MP-estimator}
  \mathcal{E}_{\dt, N}^{\mathrm{MP}2,\mathrm{BACAB}}(f) =\frac{1}{N} \sum_{n=0}^{N-1}\left(f(y^n) - \pi_{\gamma, \dt}(f)\right) z^N,
\end{equation}
with the second order weight process
{
\begin{equation}
\label{eqn:2nd-order-weight}
z^N = \sqrt{\frac{\beta}{2\gamma}} \sum_{n=0}^{N-1} 
\dt^{1/2} F(q^n)^\T G^n + \frac{1}{2}\dt^{3/2} (p^n)^\T M^{-1} \nabla_q F(q^n) G^n \,,
\end{equation}
}
obtained from $z^N=\sum_{n=0}^{N-1} \DMN(g-\dt\mathcal{L}g/2)$ by truncating 
the expansion in \eqref{eq:expansion_martingale_increment_2nd_order}. See Appendix~\ref{app:secondorder} for more detailed calculations.

We conclude this section by providing the expressions for the auxiliary discrete martingale~$z^N$ associated with other splitting schemes: 
\begin{itemize}
    \item[(i)] For the scheme $P_{\dt}^{A, B, \gamma C, B, A}$, the same expansion as~\eqref{eq:expansion_BACAB} holds, so that the second order MP estimator for this scheme coincides with~\eqref{eqn:second-order-MP-estimator}.
    \item[(ii)] For the scheme $P_{\dt}^{\gamma C, B, A, B, \gamma C}$ (recall~\eqref{eqn:CBABCB}), the counterpart of the expansion~\eqref{eq:expansion_BACAB} reads
      \begin{equation}
      \label{eq:expansion_CBABC}
      \begin{split}
        \Phi_{\dt}^q(y^n; G^n) &= \dt M^{-1} p^n + \dt^{3/2} \sqrt{\frac{\gamma}{\beta}}M^{-1} G_1^n + \Delta_{\dt}^q(y^n; G^n)\,,\\
        \Phi_{\dt}^p(y^n; G^n) 
        &= h^{1/2}\sqrt{\frac{\gamma}{\beta}} (G_1^n+G_2^n) - \dt (\nabla V(q^n) + \gamma M^{-1} p^n) \\
        & \;\;\; -\frac{\dt^{3/2}}{4}\sqrt{\frac{\gamma}{\beta}}\gamma M^{-1} (3G_1^n+G_2^n) + \Delta_{\dt}^p(y^n; G^n)\,.
\end{split}
      \end{equation}
      The expansion is in fact the same for the scheme~$P_{\dt}^{\gamma C, A, B, A, \gamma C}$. A second order MP estimator is then obtained for these two schemes by choosing (see Appendix~\ref{app:secondorder} for more detailed calculations)
     \begin{equation}\label{eqn:2nd-order-weight-CBABC}
    \begin{aligned}
    z^{n+1} - z^n & =  \frac{1}{2}\dt^{1/2}\sqrt{\frac{\beta}{\gamma}}  F(q^n)^{\T} (G_1^n + G_2^n)\\ 
    & + \frac{\dt^{3/2}}{4} \sqrt{\frac{\beta}{\gamma}} \left[ (p^n)^\T M^{-1} \nabla_q F(q^n) (G_1^n+G_2^n)+ (p^n)^\T  \nabla_q F(q^n)^\T M^{-1} (G_1^n-G_2^n) \phantom{\frac12}\right.\\ 
    & \qquad \qquad \qquad \qquad \qquad \qquad\left. - \frac{\gamma}{2} F(q^n)^\T M^{-1} (G_1^n - G_2^n)\right]\,.
    \end{aligned}
    \end{equation}
\end{itemize}
We do not consider the second order splitting schemes associated with the semigroups $P_{\dt}^{B, \gamma C, A, \gamma C, B}$ and $P_{\dt}^{A, \gamma C, B, \gamma C, A}$ since the marginal in the position variable of the invariant measures of these numerical schemes in the overdamped limit are not consistent with the Boltzmann--Gibbs measure with density proportional to~$\mathrm{e}^{-\beta V(q)}$ \cite{leimkuhler2016computation}.

\subsection{Variance analysis of martingale product estimators}
\label{subsec:MP-var-analysis}

The key advantage of the martingale product estimator is that its variance does not grow with respect to the integration time, as stated in the following result. As will be seen below, this is in sharp contrast with estimators based on Green--Kubo formulas for instance (see Theorem~\ref{thm:GK-variance}).

\begin{theorem}[Variance of MP estimators]
  \label{thm:MP-variance}
  Fix an observable $f \in \mS$, and consider any of the first ($\alpha=1$) or second ($\alpha=2$) order splitting schemes introduced in Section~\ref{sec:LD-application}, and their associated MP estimators constructed in~\eqref{eqn:first-order-MP-estimator}-\eqref{eqn:1st-order-weight} and Section~\ref{subsec:2nd-order-MP} (such as~\eqref{eq:MP2_estimator}), respectively. There exist $\dt^* > 0$ and $C \in \mathbb{R}_+$ such that, for any $\dt\in (0,\dt^*]$ and any integer $N\geq 1$ satisfying~$N\dt \geq 1$,
  \[
    \mathrm{Var}_{\dt}\left(\mathcal{E}_{\dt,N}^{\mathrm{MP}\alpha} \right) \leq C\left(1+\frac{1}{N\dt}\right)\PERIOD
\]
\end{theorem}

\begin{proof}
For the ease of notation, we write the proof for the first order MP estimator~\eqref{eqn:first-order-MP-estimator}-\eqref{eqn:1st-order-weight}, but the estimates are completely similar for second order estimators. The proof closely follows the lines of the proof of~\cite[Theorems~4.3 and~4.6]{plechac2019convergence}. Here as well, $C$ denotes a generic constant which may change from line to line. The index~$s$ of Lyapunov functions can also change from line to line.

In view of the discrete time Poisson equation~\eqref{eqn:disc-Poisson-eqn}, we can rewrite $\mathcal{E}_{\dt, N}^{\mathrm{MP}1}$ as
\[
  \mathcal{E}_{\dt, N}^{\mathrm{MP}1} = \frac{1}{N\dt} \sum_{n=0}^{N-1} \left( \WHF(y^{n+1}) - (P_{\dt} \WHF)(y^n) \right)z^N + \frac{1}{N\dt} \left( \WHF(y^0) - \WHF(y^N)  \right) z^N \COMMA
\]
so that, using $\mathrm{Var}_{\dt}\left( \mathcal{E}_{\dt, N}^{\mathrm{MP}1} \right) \leq \bE_{\dt}\left\{ \left( \mathcal{E}_{\dt, N}^{\mathrm{MP}1} \right)^2 \right\}$ and a Cauchy--Schwarz inequality,
\begin{equation}
  \label{eqn:var-estimate}
  \begin{aligned}
    \mathrm{Var}_{\dt}\left( \mathcal{E}_{\dt, N}^{\mathrm{MP}1} \right) & \leq \frac{2}{N^2\dt^2} \mathbb{E}_{\dt}\left\{\left(\sum_{n=0}^{N-1}\DMN(\WHF) z^N\right)^2  \right\} + \frac{2}{N^2\dt^2}\mathbb{E}_{\dt}\left\{ \left[\left(\WHF(y^N) - \WHF(y^0)\right)z^N \right]^2\right\} \PERIOD
  \end{aligned}
\end{equation}

We start by estimating the first term on the right hand side of~\eqref{eqn:var-estimate}. Denoting by~$\DN_n = z^{n+1} - z^n$ for notational homogeneity, a simple computation shows that 
\begin{equation*}
\begin{split}
& \bE_{\dt}\left\{\left(\sum_{n=0}^{N-1}\DMN\!\left(\WHF\right) z^N\right)^2  \right\} = \sum_{n_1, n_2, n_3, n_4=0}^{N-1}\bE_{\dt}\left\{ \DM_{n_1}\!\!\left(\WHF\right)\DM_{n_2}\!\!\left(\WHF\right)  \DN_{n_3} \DN_{n_4}\right\}\\
& = \sum_{n_1,n_2=0}^{N-1} \bE_{\dt}\left\{ \DM_{n_1}\!\!\left(\WHF\right)^2 \DN_{n_2}^2 \right\}  + 2\!\!\!\!\!\!\!\!\! \sum_{0 \leq n_1 < n_2 \leq N-1}\!\!\!\!\!\!\!\!\!\bE_{\dt}\left\{\DM_{n_1}\!\!\left(\WHF\right)\DM_{n_2}\!\!\left(\WHF\right)\DN_{n_1}\DN_{n_2}\right\}\COMMA
\end{split}
\end{equation*}
where the second equality holds because the underlying martingale increments are conditionally independent between successive steps. We next replace~$\WHF$ by~$\WTF$, up to some error term, and perform expansions in powers of~$\dt$ of $\DMN(\WTF)$. More precisely, we obtain from~\eqref{eqn:discrete-Poisson-equation-projected} and~\eqref{eqn:phi} that $-\LOPERHPI (\WTF-\WHF) = \dt^2 \phi_{1, \dt, f}$
and hence by definition of $\widetilde{\mL}_{\dt}$,
\[
\WTF-\WHF = -h^2 \mL_{\dt}^{-1} (\phi_{1, \dt, f} - \pi_{\gamma, \dt}(\phi_{1, \dt, f})).
\] 
By~\eqref{eqn:control-of-phi} and Corollary~\ref{cor:bounded-disc-Poission-solution}, there exist therefore~$\dt^*>0$, ~$C \in \mathbb{R}_+$ and~$s \in \mathbb{N}$ such that
\[
\forall \dt \in (0,\dt^*], \qquad \left|\WTF-\WHF\right| \leq C \dt^2 \mKS \PERIOD 
\]
This implies that, for any~$\dt \in (0,\dt^*]$,
\[
\left|\DMN\!\left(\WTF\right)^2-\DMN\!\left(\WHF\right)^2\right| = \left|\DMN\!\left(\WTF-\WHF\right)\right| \,  \left|\DMN\!\left(\WTF+\WHF\right)\right|\leq C \dt^2 \left[\mKS(y^n)+\mKS(y^{n+1})\right] \PERIOD
\]
so that
\begin{equation}
\label{eq:bound_difference_WTF_WHF_variance}
\begin{aligned}
& \left| \bE_{\dt}\left\{ \DM_{n_1}\!\!\left(\WHF\right)^2 \DN_{n_2}^2 \right\} - \bE_{\dt}\left\{ \DM_{n_1}\!\!\left(\WTF\right)^2 \DN_{n_2}^2 \right\} \right| \\
& \qquad \qquad \leq C\dt^3 \bE_{\dt}\left[ \left(F(q^{n_2})^{\T} G^{n_2} \right)^2\left(\mKS(y^{n_1})+\mKS(y^{n_1+1})\right)\right] \\
& \qquad \qquad \leq C\dt^3 \bE_{\dt}\left[ \left(F(q^{n_2})^{\T} G^{n_2} \right)^4 + \mKS(y^{n_1})^2 +\mKS(y^{n_1+1})^2 \right] \leq C\dt^3.
\end{aligned}
\end{equation}
We therefore obtain that, for $\dt \in (0,\dt^*]$,
\[
\left| \sum_{n_1,n_2=0}^{N-1} \bE_{\dt}\left\{ \DM_{n_1}\!\!\left(\WHF\right)^2 \DN_{n_2}^2 \right\} - \sum_{n_1,n_2=0}^{N-1} \bE_{\dt}\left\{ \DM_{n_1}\!\!\left(\WTF\right)^2 \DN_{n_2}^2 \right\} \right| \leq C N^2\dt^3.  
\]
Moreover, in view of the expansion~\eqref{eq:Taylor_WTF} of~$\DMN(\WTF)$ in powers of~$\dt$, and the bounds~\eqref{eq:moment_bounds_psi_Upsilon}, a computation similar to~\eqref{eq:bound_difference_WTF_WHF_variance} provides
\[
\frac{1}{N^2\dt^2}\sum_{n_1=0}^{N-1}\sum_{n_2=0}^{N-1} \bE_{\dt}\left\{ \DM_{n_1}\!\!\left(\WTF\right)^2 \DN_{n_2}^2 \right\} \leq C\PERIOD
\]
By summing the last two estimates, we finally obtain that, for $\dt \in (0,\dt^*]$,
\[
\frac{1}{N^2\dt^2}\sum_{n_1=0}^{N-1}\sum_{n_2=0}^{N-1} \bE_{\dt}\left\{ \DM_{n_1}\!\!\left(\WHF\right)^2 \DN_{n_2}^2 \right\} \leq C\PERIOD
\]
A similar estimate can be obtained for the other double sum 
\[
\frac{1}{N^2\dt^2}\sum_{0 \leq n_1 < n_2 \leq N-1}^{N-1}\!\!\!\!\!\!\!\!\!\bE_{\dt}\left\{\DM_{n_1}\!\!\left(\WHF\right)\DM_{n_2}\!\!\left(\WHF\right)\DN_{n_1}\DN_{n_2}\right\} \leq  C\dt\PERIOD
\]
so that the first term on the right hand side of~\eqref{eqn:var-estimate} is uniformly bounded by $C$ for $\dt$ sufficiently small.

Now, for the second term on the right hand side of~\eqref{eqn:var-estimate}, a Cauchy--Schwarz inequality gives 
\[
\mathbb{E}_{\dt}\left\{ \left[ \left(\WHF(y^N) - \WHF(y^0)\right)z^N \right]^2\right\} \leq \mathbb{E}_{\dt}\left\{ \left(\WHF(y^N) - \WHF(y^0)\right)^4 \right\}^{1/2} \mathbb{E}_{\dt}\left\{ (z^N)^4 \right\}^{1/2}.
\]
Note that by Corollary~\ref{cor:bounded-disc-Poission-solution}, $\WHF$ is uniformly bounded in~$B_{\mKS}^{\infty}(\mY)$ for some integer~$s$ and hence we can write (upon possibly increasing~$s$) 
\[
\mathbb{E}_{\dt}\left\{ \left(\WHF(y^N) - \WHF(y^0)\right)^4 \right\} \leq C \bE_{\dt}\left\{\mKS(y^0) + \mKS(y^N)\right\},
\]
which is uniformly bounded by a constant for~$\dt$ sufficiently small in view of~\eqref{eqn:K-evolution-estimate}. The estimation of $\mathbb{E}_{\dt}\left\{ (z_N)^4 \right\}$ is similar to that of the first term on the right hand side of~\eqref{eqn:var-estimate} by noting that
\[
\mathbb{E}_{\dt}\left\{ (z^N)^4 \right\} = \sum_{n=0}^{N-1} \bE_{\dt}\left\{ \DN_{n}^4 \right\} + 6\sum_{0 \leq n_1 < n_2 \leq N-1} \bE_{\dt}\left\{ \DN_{n_1}^2 \DN_{n_2}^2 \right\},
\]
which is uniformly bounded by $C N^2\dt^2$ for $\dt$ sufficiently small. Therefore, the second term on the right hand side of~\eqref{eqn:var-estimate} is uniformly bounded by $C/(N \dt)$, which concludes the desired variance estimate.
\end{proof}

\section{Comparison with the Green-Kubo method}
\label{sec:LR-vs-GK}

The Green-Kubo (GK) formula is one of the traditional methods for computing linear responses and transport coefficients (see for instance~\cite[Chapter~13]{tuckerman2010statistical}). We first recall known results on the bias of GK estimators in Section~\ref{sec:GK_bias}, and then characterize their variance in Section~\ref{sec:variance_GK}. The main conclusion is that the bias of GK estimators is comparable to the timestep bias of MP estimators, but the finite time bias of GK estimators is much smaller than the one from for MP estimators, while on the other hand the variance of GK estimators scales as~$N\dt$, in sharp contrast to the variance of MP estimators, which is uniformly bounded in time. This suggests that MP estimators could be more reliable than GK estimators.

\subsection{Green-Kubo estimators}
\label{sec:GK_bias}

Recall the reformulation~\eqref{eqn:GK} of the linear response as an integrated correlation function: 
\[
\LIN(f) = \int_0^{\infty} \bE^{\pi}\left\{ f(y_t) \varphi(y_0) \right\}\, dt = \int_0^{\infty} \bE^{\pi}\left\{ \left[f(y_t)-\pi(f)\right] \varphi(y_0) \right\}\, dt,
\]
where 
\begin{equation}
\label{eq:def_varphi_GK}
\varphi = \LOPERTIL^* \mathbf{1}
\end{equation}
has average~0 with respect to~$\pi$ and~$\bE^{\pi}$ is the expectation taken with respect to all realizations of the Markov process $\{y_t\}_{t \geq 0}$ and with respect to all initial conditions~$y_0$ distributed according to~$\pi$. The timestep bias of GK estimators was made precise in~\cite{leimkuhler2016computation,lelievre2016partial}. The result essentially states that the linear response can be approximated by a Riemann sum, or using a trapezoidal rule for second order schemes. To make this result precise, we introduce the following estimators:
\[
\mathcal{E}_{\dt, N}^{\mathrm{GK}1}(f) = \dt \sum_{n=0}^{N-1} \left[f(y^n) - \pi_{\gamma, \dt}(f)\right] \varphi(y^0), \qquad y^0 \sim \pi_{\gamma, \dt}\COMMA
\]
and
\[
\mathcal{E}_{\dt, N}^{\mathrm{GK}2}(f) = \dt \left( \frac12 \left[f(y^0) - \pi_{\gamma, \dt}(f)\right] + \sum_{n=1}^{N-1} \left[f(y^n) - \pi_{\gamma, \dt}(f)\right] \right)\varphi(y^0), \qquad y^0 \sim \pi_{\gamma, \dt} \PERIOD
\]
Note that, as discussed after~\eqref{eqn:desired-estimator}, we consider here as well a ``perfect'' centering based on subtracting~$\pi_{\gamma, \dt}(f)$; but the same discussion as for MP estimators applies here. Let us also mention that the initial distribution for~$y^0$ has to be the invariant probability measure~$\pi_{\gamma, \dt}$ for GK estimators whereas MP estimators do not require this condition. 

We can now recall the results from~\cite{leimkuhler2016computation,lelievre2016partial}, which correspond to the case $N \to \infty$. We denote by $\bE_{\dt}^{\pi_{\gamma, \dt}}$ the expectation taken with respect to all realizations of the Markov chain $\{y^n\}_{n \geq 0}$ and with respect to all initial conditions~$y^0$ distributed according to~$\pi_{\gamma,\dt}$. 

\begin{theorem}
  \label{thm:error-estimates-GK}
  Fix an observable $f \in \mS$, and consider any of the first ($\alpha=1$) or second ($\alpha=2$) order splitting schemes introduced in Section~\ref{sec:LD-application}. Then there exist~$\dt^* > 0$ and $C \in \mathbb{R}_+$ such that, for any $\dt\in(0,\dt^*]$,
\[
  \left| \bE_{\dt}^{\pi_{\gamma, \dt}}\left\{\mathcal{E}_{\dt, \infty}^{\mathrm{GK}\alpha}(f)\right\} - \LIN(f) \right| \leq C \dt^\alpha. 
\]
\end{theorem}

In practice, the integrated correlation over the infinite time horizon has to be approximated by that over a finite time horizon, which adds another bias, exponentially small in the integration time~$N\dt$ -- and hence asymptotically much smaller than the finite time integration bias of order~$1/\sqrt{N\dt}$ for MP estimators.

\begin{corollary}[Bias of GK estimators]
  \label{cor:GK-bias}
  Fix an observable $f \in \mS$, and consider any of the first ($\alpha=1$) or second ($\alpha=2$) order splitting schemes introduced in Section~\ref{sec:LD-application}. Then there exist $\lambda,\dt^* > 0$ and $C \in \mathbb{R}_+$ such that, for any $\dt\in(0,\dt^*]$, 
\[
  \forall N \geq 1, \qquad \left| \bE_{\dt}^{\pi_{\gamma, \dt}}\left\{\mathcal{E}_{\dt,N}^{\mathrm{GK}\alpha}(f)\right\} - \LIN(f) \right| \leq C \left( \dt^\alpha + \mathrm{e}^{-\lambda N\dt} \right) \PERIOD 
\]
\end{corollary}

\begin{proof}
  We bound the error as
  \[
  \begin{aligned}
    \left| \bE_{\dt}^{\pi_{\gamma, \dt}}\left\{\mathcal{E}_{\dt,N}^{\mathrm{GK}\alpha}(f)\right\} - \LIN(f) \right| & \leq \left| \bE_{\dt}^{\pi_{\gamma, \dt}}\left\{\mathcal{E}_{\dt,N}^{\mathrm{GK}\alpha}(f)\right\} - \bE_{\dt}^{\pi_{\gamma, \dt}}\left\{\mathcal{E}_{\dt,\infty}^{\mathrm{GK}\alpha}(f)\right\} \right| \\
    & \quad + \left| \bE_{\dt}^{\pi_{\gamma, \dt}}\left\{\mathcal{E}_{\dt,\infty}^{\mathrm{GK}\alpha}(f)\right\} - \LIN(f) \right|.
    \end{aligned}
  \]
  The second term on the right hand side of the previous inequality is bounded by~$C\dt^{\alpha}$ in view of Theorem~\ref{thm:error-estimates-GK}. Choosing an integer~$s$ such that $f,\varphi \in B^\infty_{\mKS}(\mY)$, and using the exponential decay estimates of Theorem~\ref{thm:ergodicity-splitting}, we obtain, for some constant~$K \in \mathbb{R}_+$,
\[
  \begin{split}
    \left| \bE_{\dt}^{\pi_{\gamma, \dt}}\left\{\mathcal{E}_{\dt,N}^{\mathrm{GK}\alpha}(f)\right\} - \bE_{\dt}^{\pi_{\gamma, \dt}}\left\{\mathcal{E}_{\dt,\infty}^{\mathrm{GK}\alpha}(f)\right\} \right| & \leq \dt \sum_{n = N}^{\infty}\bE_{\dt}^{\pi_{\gamma, \dt}} \left( \left| \left[f(y^n) - \pi_{\gamma, \dt}(f)\right] \varphi(y^0) \right| \right)\\
    & \leq K \|f\|_{B^\infty_{\mKS}} \|\varphi\|_{B^\infty_{\mKS}} \pi_{\gamma, \dt}\left(\mKS^2\right) \dt \sum_{n = N}^{\infty} \mathrm{e}^{-\lambda n \dt} 
  \end{split}
\]
The result follows by noting that $\dt\sum_{n = N}^{\infty}\mathrm{e}^{-\lambda n \dt} \leq 2\lambda^{-1} \mathrm{e}^{-\lambda N \dt}$ for~$\dt$ sufficiently small.  
\end{proof}

Let us conclude this section by emphasizing that all the results presented here require that initial conditions are distributed according to the invariant probability measure~$\pi_{\gamma,\dt}$, in contrast to the estimators constructed in Section~\ref{sec:general-continuous-setting}.

\subsection{Variance analysis of Green-Kubo estimators}
\label{sec:variance_GK}

We can now present the main result of this section, which states that the variance of the Green--Kubo estimator is given at dominant order by the variance of the random variable
\[
\mathscr{E}_{\dt, N}^{\mathrm{GK}\alpha}(f) = \varphi(y^0) \sum_{n=0}^{N-1} \DMN\left(\WHF \right) \COMMA 
\]
which is expected to scale linearly with the integration time~$Nh$. Recall that~$\varphi$ is given by~\eqref{eq:def_varphi_GK}. 
The increase of the variance with time can significantly deteriorate the performance of the estimator especially when correlation functions decay slowly.

\begin{theorem}[Variance of GK estimators]
  \label{thm:GK-variance}
  Fix an observable $f \in \mS$, and consider any of the first ($\alpha=1$) or second ($\alpha=2$) order splitting schemes introduced in Section~\ref{sec:LD-application}. Then there exist $\dt^* > 0$ and $C \in \mathbb{R}_+$ such that, for any $\dt\in(0,\dt^*]$,
    \[
        \forall N \geq 1, \qquad \left| \mathrm{Var}_{\dt}^{\pi_{\gamma, \dt}}\left(  \mathcal{E}_{\dt, N}^{\mathrm{GK}\alpha}(f)  \right)- 
        \mathrm{Var}_{\dt}^{\pi_{\gamma, \dt}}\left( \mathscr{E}_{\dt, N}^{\mathrm{GK}\alpha}(f) \right)
        \right| \leq C\left(1+ \sqrt{N\dt}\right).
    \]
    As a consequence, there exists a constant~$K \in \mathbb{R}_+$ such that, for any $\dt\in(0,\dt^*]$, the following bound holds when~$N\dt \geq 1$:
    \[
    \mathrm{Var}_{\dt}^{\pi_{\gamma, \dt}}\left(  \mathcal{E}_{\dt, N}^{\mathrm{GK}\alpha}(f)\right) \leq K N\dt.
    \]
\end{theorem}

\begin{proof}
  The proof is quite similar to the proof of Theorem~\ref{thm:MP-variance}, which is not surprising since the GK estimator shares some structural similarities with MP estimators (in fact, replacing~$z^N/(N\dt)$ by~$\varphi(y^0)$). We therefore follow the notation introduced in Theorem~\ref{thm:MP-variance}. Note first that
  \[
  \mathcal{E}_{\dt, N}^{\mathrm{GK}\alpha}(f) = \mathscr{E}_{\dt, N}^{\mathrm{GK}\alpha}(f) + \mathscr{R}_{\dt, N}^{\mathrm{GK}\alpha}(f), \qquad \mathscr{R}_{\dt, N}^{\mathrm{GK}\alpha}(f) = \left( \WHF(y^0) - \WHF(y^N)  \right) \varphi(y^0) \PERIOD 
  \]
  Therefore, in view of the inequality $|\mathrm{Var}(X+Y)-\mathrm{Var}(X)| \leq \mathrm{Var}(Y) + 2\sqrt{\mathrm{Var}(Y)\mathrm{Var}(X)}$, we obtain
  \[
  \begin{aligned}
    \left| \mathrm{Var}_{\dt}^{\pi_{\gamma, \dt}}\left(  \mathcal{E}_{\dt, N}^{\mathrm{GK}\alpha}(f)  \right) - \mathrm{Var}_{\dt}^{\pi_{\gamma, \dt}}\left(  \mathscr{E}_{\dt, N}^{\mathrm{GK}\alpha}(f) \right) \right| & \leq \mathrm{Var}_{\dt}^{\pi_{\gamma, \dt}}\left(\mathscr{R}_{\dt, N}^{\mathrm{GK}\alpha}(f)\right) \\ & \quad + 2 \sqrt{\mathrm{Var}_{\dt}^{\pi_{\gamma, \dt}}\left(\mathscr{R}_{\dt, N}^{\mathrm{GK}\alpha}(f)\right)\mathrm{Var}_{\dt}^{\pi_{\gamma, \dt}}\left(  \mathscr{E}_{\dt, N}^{\mathrm{GK}\alpha}(f) \right)}.
  \end{aligned}
  \]
    The estimate of the variance of~$\mathscr{E}_{\dt, N}^{\mathrm{GK}\alpha}(f)$ is similar to that of $\sum_{n=0}^{N-1}\DMN(\WHF) z^N$ in Theorem~\ref{thm:MP-variance} and can be shown to be bounded by $CN\dt$.
  Furthermore, a simple computation shows that
  \[
  \mathrm{Var}_{\dt}^{\pi_{\gamma, \dt}}\left(\mathscr{R}_{\dt, N}^{\mathrm{GK}\alpha}(f)\right)  \leq \mathbb{E}_{\dt}^{\pi_{\gamma, \dt}}\left(\mathscr{R}_{\dt, N}^{\mathrm{GK}\alpha}(f)^2\right) \leq 2 \int_\mY \left[ \WHF^2 + \left(P_\dt^N \WHF\right)^2 \right] \varphi^2 \, d\pi_{\gamma, \dt}\PERIOD
  \]
  The right hand side of the latter inequality is uniformly bounded in~$N$ in view of~\eqref{eqn:K-evolution-estimate} since~$\varphi,\WHF$ are both in~$B^\infty_{\mKS}(\mY)$ for some integer~$s \in \mathbb{N}$ (with uniform bounds for~$\WHF$ by Corollary~\ref{cor:bounded-disc-Poission-solution}). 
\end{proof}

\section{Computational benchmarks}
\label{sec:numerics}

We demonstrate the accuracy and efficiency of the derived MP estimators
by applying them to  simple two dimensional Langevin systems. 
Denoting by~$q =(q_1, q_2)$ and $p = (p_1, p_2)$, we consider 
the Langevin dynamics~\eqref{eqn:ULD} with $\beta = \gamma = 1$ and $M = \mathrm{Id}$.
We study the behavior of the proposed MP estimators in the case of unbounded position space and two different types of perturbing forces (Examples~1 and~2). In Example~3, we present comparison 
between the MP and Green-Kubo estimators in the case of periodic (compact) position domain. 

\ifx\IMAJNA\undefined
\subsubsection*{Example 1: position dependent perturbation $F(q)$.}
\else
\medskip\noindent{\it Example 1:} {position dependent perturbation $F(q)$.}
\fi
In this example we choose a harmonic potential for $q \in\mathbb{R}^2$:
\[
V(q,\VPAR) = \frac{1}{2}\VPAR |q|^2 = \frac{1}{2}\VPAR (q_1^2 + q_2^2),
\]
with~$\VPAR=1$ in our simulations. Note that the dynamics is defined on $\mathcal{Y}=\R^2\times \R^2$, thus strictly speaking not satisfying our set of assumptions. Recall however that the restriction to a compact position domain is for the ease of presentation, as hinted at after Assumption~\ref{ass:V_smooth_and_domain_compact}.

\begin{figure}[ht]
    \includegraphics[width=0.95\textwidth]{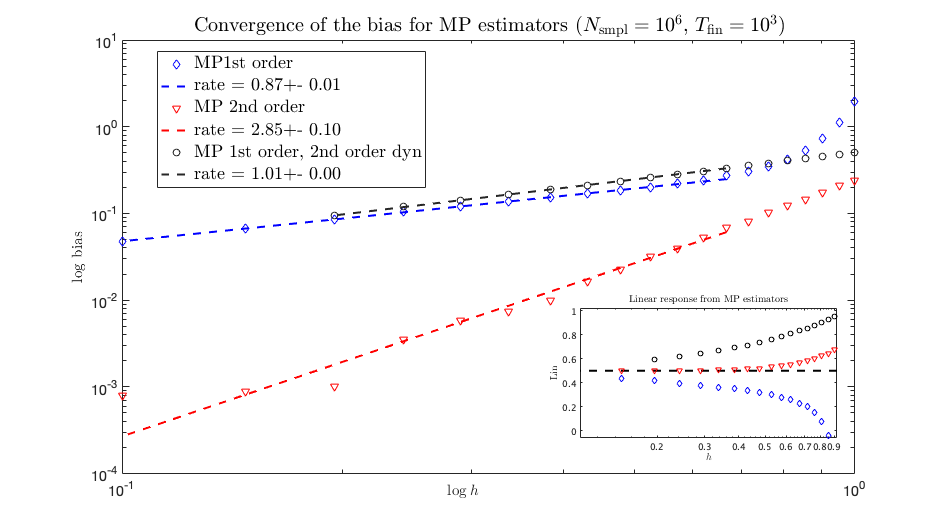}
        \caption{{\it Example 1.} Convergence rates of the bias for the MP estimators, $\mathcal{E}_{\dt, N}^{\mathrm{MP}\alpha}$, illustrating the error estimates~\eqref{eqn:1st-order-MP-bias} and~\eqref{eqn:2nd-order-MP-bias}. The inset depicts the dependence of the bias on the timestep $\dt$. The exact value $\LIN(f)=0.5$ of the sensitivity index is computed from Gaussian integrals.}
        \label{fig:MPquadconv}
\end{figure}

The perturbing force $F$ and the observable $f$ are chosen as
\[
   F_i(q)\equiv\frac{\partial^2 V}{\partial q_i\partial\VPAR} (q,\VPAR)
   = q_i\,\COMMA \qquad f(q) = |q|^2\,.
 \]
 Thus we estimate the sensitivity of the second moment of the $q$-marginal for the invariant distribution (a Gaussian distribution in $\R^2$) with respect to the parameter $\VPAR$. 
 Since the perturbing force is only position dependent, one can apply directly the formula
 \eqref{eqn:2nd-order-weight}, \emph{i.e.}, the increment of the re-weighting process is
 \begin{equation} \label{eqn:2nd-order-weight-inc}
z^{n+1}-z^n = \sqrt{\frac{\beta}{2\gamma}} \Big[ \dt^{1/2} F(q^n)^{\T} + \frac{1}{2}\dt^{3/2} (p^n)^\T M^{-1} \nabla_q F(q^n)\Big] G^n \,.
\end{equation}
 The example allows us to demonstrate, in a simple setting, the importance of the $\dt^{3/2}$ term in the definition of re-weighting process $z^n$ in \eqref{eqn:2nd-order-weight}. Simulations are run with~$\NSMPL = 10^6$ realizations, for dynamics integrated up to a physical time~$\TFIN = 10^3$. Figure~\ref{fig:MPquadconv} depicts the convergence of the bias for the sensitivity index $\LIN{f}$ for the 1st order ($\alpha=1$) and 2nd order ($\alpha=2$) MP estimators~$\mathcal{E}_{\dt, N}^{\mathrm{MP}\alpha}$, confirming the error estimates~\eqref{eqn:1st-order-MP-bias} and~\eqref{eqn:2nd-order-MP-bias}. The 2nd order estimator uses the BACAB splitting~\eqref{eqn:BACAB} for the Langevin dynamics. In contrast, an MP estimator (here called 1st order with the 2nd order dynamics) omitting the $\dt^{3/2}$ term in \eqref{eqn:2nd-order-weight-inc} 
achieves only a first order rate of convergence.

\ifx\IMAJNA\undefined
\subsubsection*{Example 2: momentum dependent perturbation $F(p)$.}
\else
\medskip\noindent{\it Example 2:} {momentum dependent perturbation $F(p)$.}
\fi
The purpose of this benchmark is to demonstrate the full expression for the 2nd-order increment of the re-weighting process in the MP estimator 
$\mathcal{E}_{\dt, N}^{\mathrm{MP}\alpha}$. We consider the same potential as in Example 1, 
and denote by%
\[
\pi_\beta(dq\,dp) = \frac{\beta^2 \omega}{4\pi^2} \mathrm{e}^{-\beta |p|^2/2} \mathrm{e}^{-\beta\VPAR |q|^2/2}\,dq\,dp\,
\]
the invariant probability measure of the Langevin dynamics~\eqref{eqn:ULD}. We modify the perturbing force as $F(p) = p$, which corresponds to a perturbation of the friction coefficient $\gamma$ in the Langevin equation~\eqref{eqn:ULD}. In view of the fluctuation-dissipation relation $\sigma\sigma^T = 2\beta^{-1} \gamma$, which relates the
diffusion matrix $S=\sigma\sigma^T$ with the friction matrix $\gamma$ and the inverse temperature~$\beta$, the perturbation of $\gamma$ results in a perturbation of the average~$\pi_\beta(f)$ with respect to $\beta$. Therefore, for a suitable choice of observables, \emph{e.g.}, $f_1(q) = q_1^2 + q_2^2$ and $f_2(p)=p_1^4 + p_2^4$ used in this example, the sensitivity/linear response $d\pi_\beta(f)/d\beta$ can be computed exactly by evaluating Gaussian integrals.

\begin{figure}[ht]
    \includegraphics[width=0.95\textwidth]{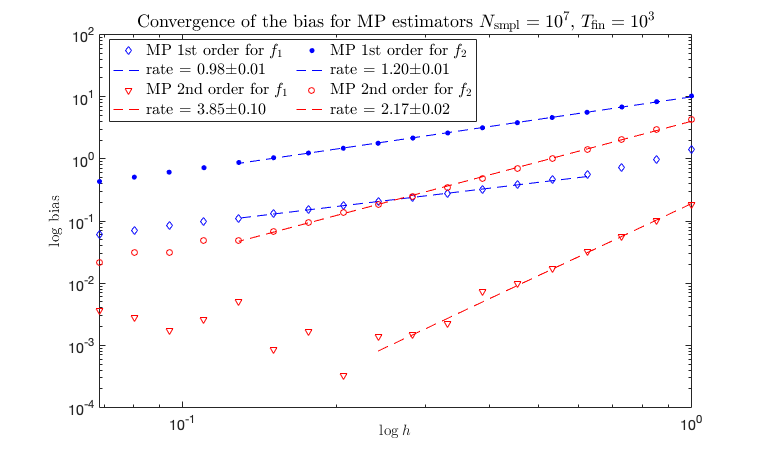}
        \caption{{\it Example 2.} Convergence rates of the bias for the MP estimators, $\mathcal{E}_{\dt, N}^{\mathrm{MP}\alpha}$. 
        The exact values $\LIN(f_1)=-2$ and $\LIN(f_2)=-12$ of the sensitivity index with respect to $\beta$ are computed from Gaussian integrals.}
        \label{fig:MPbetaconv}
\end{figure}

A simple computation based on~\eqref{eqn:gfun} shows that~$g(q,p) = \tfrac{\beta}{2\gamma}|p|^2$,
which allows us to obtain an explicit expression 
for the 
martingale increment of the 2nd-order 
MP estimator by substituting the formula for~$g$ into~\eqref{eqn:zinc} in Appendix~\ref{app:secondorder}. The results presented in Figure~\ref{fig:MPbetaconv}, obtained with~$\NSMPL = 10^7$ realizations of the dynamics integrated up to a physical time~$\TFIN = 10^3$, show that the second order MP estimators allow for passing from a bias of order~$\dt$ to a bias of order~$\dt^2$, and even to~$\dt^4$ for the observable~$f_1$. For the latter case, the bias cannot be decreased below the value~$10^{-3}$ since the statistical error starts to dominate the overall error.

\begin{figure}[ht]
    \includegraphics[width=0.95\textwidth]{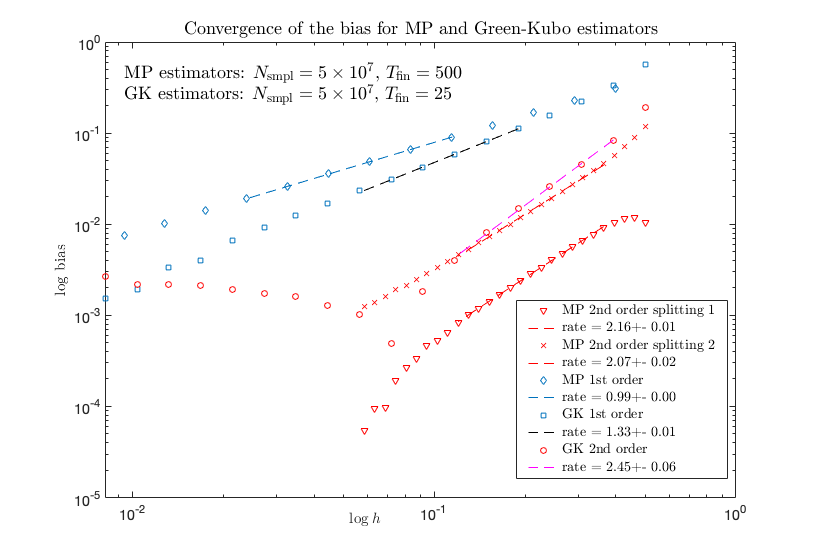}
        \caption{{\it Example 3.} Convergence rates of the bias for the different estimators. The 2nd order MP estimator is tested with the splitting schemes BACAB (splitting 1) and CBABC (splitting 2)  }\label{fig:GKMPconv}
\end{figure}

\ifx\IMAJNA\undefined
\subsubsection*{Example 3.}
\else
\medskip\noindent{\it Example 3.}
\fi
In the last example we consider a Langevin dynamics on the phase space $\mathcal{Y} = \mathbb{T}^2 \times \R^2$ with a non-conservative constant perturbing force $F \in \R^2$ (which does not derive from the gradient of a smooth periodic potential). More precisely, using the Langevin dynamics \eqref{eqn:ULD}, we sample the invariant distribution for 
$q =(q_1, q_2) \in \mathcal{D} = (2\pi \mathbb{T})^2$ and $p = (p_1, p_2) \in \mathbb{R}^2$, 
with the periodic potential 
\[
V(q) = 2 \cos(2q_1) + \cos(q_2)\,.
\] 
We are interested in estimating the mobility in a direction given by the constant unit vector $F$.  The observable $f$ and the corresponding function $\varphi$ in~\eqref{eqn:GK} are
\[
f(q, p) = F^{\T} M^{-1} p\,, \;\;\;\;\; 
\varphi(q, p) = \beta F^{\T} M^{-1} p\,.
\]
The Green-Kubo formula for the mobility therefore reduces to the following time integral of the velocity autocorrelation function: 
\[
\LIN(f)= \beta \int_{0}^{\infty} \bE^{\pi}\left\{\left(F^{\T} M^{-1} p_t\right) \left(F^{\T} M^{-1} p_0 \right)\right\} dt\PERIOD
\]
For the simulation results presented here, the perturbation is chosen to be $F = (1, 0)^{\T}$, which corresponds to estimating the mobility in the $q_1$-direction. Figure~\ref{fig:GKMPconv} shows that the predicted convergence rates for the bias of both martingale product ($\mathcal{E}_{\dt, N}^{\mathrm{MP}\alpha}$) and Green-Kubo ($\mathcal{E}_{\dt, N}^{\mathrm{GK}\alpha}$) estimators are indeed of the order predicted by our theoretical analysis. These results were produced with~$\NSMPL = 10^7$ realizations of the dynamics integrated up to a physical time~$\TFIN = 900$ for~MP estimators, and~$\TFIN = 25$ for GK estimators.

\begin{figure}[ht]
        \centering
         \includegraphics[width=0.80\textwidth]{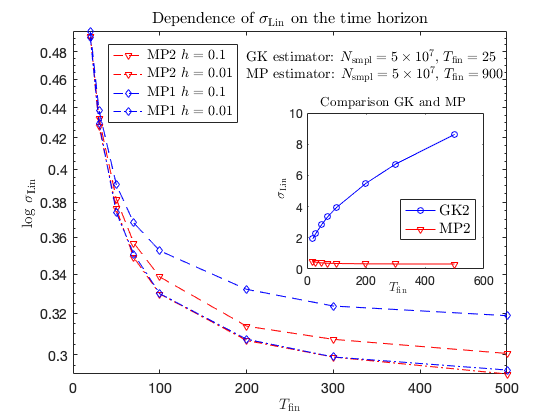}
    \caption{{\it Example 3.} Comparison variance for the martingale product (MP) and Green-Kubo (GK) estimators.}
    \label{fig:GKMPvar}
\end{figure}

An important property of the proposed MP estimators is their bounded variance while the variance of the GK estimators grows in time as demonstrated in the inset of Figure~\ref{fig:GKMPvar}. However, a fair comparison of efficiency should take into account the fact that auto-correlation functions can be decaying fast (as in the example studied here) and thus the GK estimator may not require long trajectories to be integrated, see Corollary~\ref{cor:GK-bias} for the precise statement.
On the other hand, initial conditions for the Green-Kubo estimator need
to be sampled from the 
invariant distribution and thus preparation of an initial state may still require a long time integration of the Langevin dynamics or the use of another equilibration algorithm, thus contributing to the overall computational cost of the method.  

\appendix
\section{Some useful estimates}
\label{app:useful_estimates}

We present some technical estimates which are crucial for justifying the consistency of MP estimators. Similar estimates were obtained in~\cite{plechac2019convergence} for overdamped Langevin dynamics with a compact state space setting. Throughout this section, we consider a deterministic initial condition~$y$ and denote the expectation with respect to all realizations of the Markov chain $\{y^n\}_{n \geq 1}$ starting from~$y^0$ by~$\bE_{\dt}^{y^0}$.

\begin{lemma}\label{lem:second-moment-bound}
Consider any of the splitting schemes introduced in Section~\ref{sec:splitting_Lang}, and fix an integer~$s \geq 1$. There exists constants $C \in \mathbb{R}_+$ and $\dt^* > 0$ such that, for any $0 < \dt \leq \dt^*$ and any observable $f \in B_{\mKS}^{\infty}(\mY)$, it holds
\begin{equation}
  \label{eq:lemma_A1}
  \forall n \in \mathbb{N}, \qquad \bE_{\dt}^y \left\{ \left|\WHF(y^n)\right|^2 \right\} + \bE_{\dt}^y \left\{ \left|P_\dt\WHF(y^n)\right|^2 \right\} \leq C \mKS^2(y)\left\|f\right\|_{B_{\mKS}^{\infty}}^2\COMMA
\end{equation}
where $\WHF$ is the solution to the discrete Poisson equation~\eqref{eqn:disc-Poisson-eqn}.
\end{lemma}

\begin{proof}
  In this proof and the following ones, we denote by $C$ a generic constant that can vary line by line, and fix $f \in B_{\mKS}^{\infty}(\mY)$. By Corollary~\ref{cor:bounded-disc-Poission-solution}, there exists~$\dt^*>0$ such that, for all $0 < \dt \leq \dt^*$,
  \begin{equation}
    \label{eqn:second-moment-estimate}
    \left\| |\WHF|^2 \right \|_{B_{\mKS^2}^{\infty}} = \left\|\WHF\right\|_{B_{\mKS}^{\infty}}^2 \leq \left(\frac{2K}{\lambda}\right)^2 \|f - \pi_{\gamma, \dt}(f)\|_{B_{\mKS}^{\infty}}^2.
  \end{equation}
  Therefore, 
  \[
    \bE_{\dt}^y \left\{\left|\WHF(y^n)\right|^2 \right\} \leq \left(\frac{2K}{\lambda}\right)^2 \|f - \pi_{\gamma, \dt}(f)\|_{B_{\mKS}^{\infty}}^2 \bE_{\dt}^y\{\mKS^2(y^n)\} \leq C \|f - \pi_{\gamma, \dt}(f)\|_{B_{\mKS}^{\infty}}^2 \mKS^2(y) \COMMA
  \]
  where the last inequality follows from~\eqref{eqn:K-evolution-estimate}. To bound the first term on the left hand side of~\eqref{eq:lemma_A1}, it remains to note that $\|f - \pi_{\gamma, \dt}(f)\|_{B_{\mKS}^{\infty}} \leq C \|f\|_{B_{\mKS}^{\infty}}$ in view of~\eqref{eqn:integrate-scale-function} and the inequality $|\pi_{\gamma, \dt}(f)| \leq \|f\|_{L_{\mKS}^{\infty}} \pi_{\gamma, \dt}(\mKS)$.

  For the second term on the left hand side of~\eqref{eq:lemma_A1}, the proof follows the same lines since
  \[
    \left|P_\dt\widehat{f}_{\dt}(y^n)\right| \leq \left\|\widehat{f}\right\|_{B^\infty_{\mKS}} (P_\dt \mKS)(y^n) \leq C \left\|\widehat{f}\right\|_{B^\infty_{\mKS}} \mKS(y^n) \COMMA
  \]
  where the last inequality is obtained from~\eqref{eqn:K-evolution-estimate}. 
\end{proof}

\begin{proposition}
  \label{prop:elementary-term}
  Consider any of the splitting schemes introduced in Section~\ref{sec:splitting_Lang}, and fix an integer~$s \geq 1$. Assume that $\{z^n\}_{n \geq 0}$ (which depends on the realization~$\{y^n\}_{n \geq 0}$ of the numerical scheme) is a zero mean martingale with independent increments~$\{z^{n+1}-z^n\}_{n \geq 0}$, such that~$z^0 = 0$, and which satisfies the following bound: There exist~$D \in \mathbb{R}_+$ and~$\dt^*>0$ such that
  \begin{equation}
  \label{eq:ass_eta_increment}
  \forall n \in \mathbb{N}, \qquad \bE_{\dt}^y \left\{\left(z^{n+1}  - z^n\right)^2 \right\} \leq D\dt \PERIOD
  \end{equation}
  Upon possibly reducing~$\dt^*$, there exists $C \in \mathbb{R}_+$ such that, for any~$\dt \in (0,\dt^*]$ and~$\phi \in B_{\mKS}^{\infty}(\mY)$,
  \begin{equation}
    \label{eqn:crude-estimate-1}
    \left|\frac{1}{N}\sum_{n=0}^{N-1}\bE_{\dt}^y\left\{ \phi(y^n)z^N\right\}\right| \leq \frac{C}{\sqrt{\dt}}\mKS(y) \left\|\phi\right\|_{B_{\mKS}^{\infty}} \PERIOD
\end{equation}
\end{proposition}

The interest of the bound~\eqref{eqn:crude-estimate-1}, which is crude in terms of the dependence on~$h$, is that it explicitly depends on the norm of~$\phi$, and can therefore be applied to families of functions indexed by~$\dt$.

\begin{proof}
From the continuous time Poisson equation~\eqref{eqn:cont-Poisson-eqn}, upon rearranging terms as in~\eqref{eqn:1st-order-estimate} and recalling~\eqref{eqn:obs_increment}, 
\[
\begin{aligned}
  \frac{1}{N}\sum_{n=0}^{N-1}\bE_{\dt}^y\left\{ \phi(y^n)z^N\right\} & = \frac{1}{N\dt}\sum_{n=0}^{N-1}\bE_{\dt}^y \left\{ \DMN(\widehat{\phi}_{\dt}) z^N\right\} + \frac{1}{N\dt}\bE_{\dt}^y \left\{ \left(\widehat{\phi}_{\dt}(y) - \widehat{\phi}_{\dt}(y^N)\right) z^N \right\} \\
  & \!\!= \frac{1}{N\dt}\sum_{n=0}^{N-1}\bE_{\dt}^y \left\{ \DMN(\widehat{\phi}_{\dt}) \left(z^{n+1}-z^n\right)\right\} + \frac{1}{N\dt}\bE_{\dt}^y \left\{ \left(\widehat{\phi}_{\dt}(y) - \widehat{\phi}_{\dt}(y^N)\right) z^N \right\},
\end{aligned}
\]
so that, by the Cauchy--Schwarz inequality,
\[
  \begin{aligned}
    \frac{1}{N} \left| \sum_{n=0}^{N-1} \bE_{\dt}^y\left\{ \phi(y^n)z^N\right\} \right| & \leq \sqrt{\frac{D}{N\dt}} \left(\sum_{n=0}^{N-1}\bE_{\dt}^y \left\{ \left[\DMN(\widehat{\phi}_{\dt})\right]^2 \right\} \right)^{1/2} \\
    & \ \ \ \ + \frac{2}{N\dt}\left( \bE_{\dt}^y\left\{\left(z^N\right)^2\right\} \right)^{1/2} \left( \left|\widehat{\phi}_{\dt}(y)\right|^2  + \bE_{\dt}^y \left\{\left|\widehat{\phi}_{\dt}(y^N)\right|^2\right\} \right)^{1/2} \PERIOD
  \end{aligned}
\]
The second term on the right hand side of the previous inequality is bounded using Lemma~\ref{lem:second-moment-bound} and~\eqref{eq:ass_eta_increment}, which implies that
\[
\bE_{\dt}^y\left\{\left(z^N\right)^2\right\} = \sum_{n=0}^{N-1} \bE_{\dt}^y \left\{\left(z^{n+1}  - z^n\right)^2 \right\} \leq DN\dt ;
\]
while we use a Cauchy--Schwarz inequality and Lemma~\ref{lem:second-moment-bound} again for the first one:
\[
  \bE_{\dt}^y \left\{ \left[\DMN(\widehat{\phi}_{\dt})\right]^2 \right\} \leq 2 \left( \bE_{\dt}^y \left\{ \left[\widehat{\phi}_{\dt}(y^n)\right]^2 \right\} + \bE_{\dt}^y \left\{ \left[P_\dt\widehat{\phi}_{\dt}(y^n)\right]^2 \right\} \right) \leq C \mKS^2(y)\left\|\phi\right\|_{B_{\mKS}^{\infty}}^2 \PERIOD
\]
The desired estimate finally follows by combining the above two estimates.
\end{proof}

\begin{remark}
{\rm
For overdamped Langevin dynamics with compact state space, where $\mKS \equiv 1$ for all $s \geq 1$, the estimate~\eqref{eqn:crude-estimate-1} reduces to
\[
\left|\frac{1}{N}\sum_{n=0}^{N-1}\bE_{\dt}\left\{ \phi(y^n)z^N\right\}\right| \leq \frac{C}{\sqrt{\dt}}\left\|\phi\right\|_{B^{\infty}},
\]
which is consistent with the estimate obtained in~\cite[Lemma~6.1]{plechac2019convergence}.
}
\end{remark}

\section{Discrete re-weighting martingale for second-order schemes.}
\label{app:secondorder}

Following the general methodology of Section~\ref{sec:construction_discrete_martingale}
we first derive in this section the auxiliary discrete martingale~$z^n$ associated with the second order BACAB splitting scheme in the case when $F$ depends on both $q$ and $p$. 
Expressions for the second order CBABC splitting are given at the end of this section but only for the case when $F$ depends on $q$. 

The calculation requires to expand $\DMN(g - \dt \mL g/2)$ up to order ${h}^{3/2}$ in~\eqref{eq:MP2_estimator}. 
We use to this end the expansion~\eqref{eq:expansion_BACAB} of the one-step increment of the BACAB scheme up to the order~$\dt^{3/2}$.
An Ito--Taylor formula allows us to expand $g(y^{n+1}) = g(y^n + \Phi_{\dt}(y^n; G^n))$ as
\begin{equation}\label{eq:series}
\begin{split}
g(y^{n+1}) = g(y^n) & + \nabla g(y^n)^{\T} \Phi_{\dt}(y^n; G^n) + \frac{1}{2} D^2 g(y^n) : \Phi_{\dt}(y^n; G^n)^{\otimes 2} \\
           & +\frac{1}{6} D^3 g(y^n) : \Phi_{\dt}(y^n; G^n)^{\otimes 3} + \ldots \COMMA
\end{split}
\end{equation}
where (for $y=(q,p)\in\R^{D}\times \R^D$)
\[
D^k g(y):v_1 \otimes \dots \otimes v_k \equiv \sum_{i_1,\dots,i_k=1}^{2D} \frac{\partial^k g(y)}{\partial y_{i_1} \dots \partial y_{i_k}} v_{1,i_1}\dots v_{k,i_k} \,.
\]
{We also recall that, with our convention, $\nabla F$ is a matrix with entries
$(\nabla F)_{ij} = \partial_{y_i} F_j$, i.e., the columns of~$\nabla F$ are
the gradients of the components of the vector field~$F$. Furthermore, the gradient of a scalar valued function is understood as a column vector.}
Combining~\eqref{eq:series} with the expansion~\eqref{eq:expansion_P_dt} of $P_{\dt}$, we next expand $g(y^{n+1}) - (P_{\dt}g)(y^n)$ up to terms of order~$\dt^{3/2}$, and collect terms of corresponding powers as follows: 
\begin{equation}\label{eqn:part-1}
\begin{split}
\mathcal{O}(\dt^{1/2}): & \;  \sqrt{\frac{2\gamma}{\beta}}  \nabla_p g(y^n)^{\T} G^n \,\COMMA
\\
\mathcal{O}(\dt): & \; \frac{\gamma}{\beta} D_{pp}^2 g(y^n) : (G^n)^{\otimes 2} 
- \frac{\gamma}{\beta} \Delta_p g(y^n) \COMMA\\
\mathcal{O}(\dt^{3/2}): & \sqrt{\frac{\gamma}{2\beta}}\nabla_q g(y^n)^{\T} M^{-1}G^n -\sqrt{\frac{\gamma}{2\beta}} \gamma \nabla_p g(y^n)^{\T}  M^{-1} G^n \\
& \quad +\sqrt{\frac{2\gamma}{\beta}} D_{qp}^2 g(y^n) : (M^{-1} p^n) \otimes G^n -\sqrt{\frac{2\gamma}{\beta}} D_{pp}^2 g(y^n) : (\nabla V(q^n) + \gamma M^{-1} p^n) \otimes G^n\\
& \quad + \frac{1}{6} \left(\frac{2\gamma}{\beta}\right)^{3/2} D_{ppp}^3 g(y^n) : (G^n)^{\otimes 3} \PERIOD
\end{split}
\end{equation}
We obtain the expansion of $\DMN(g - \dt\mL g/2)$ by replacing $g$ with $g-\dt\mL g/2$ and truncating the resulting expression at order~$\dt^{3/2}$: 
\begin{equation}\label{eq:incrementZ}
z^{n+1} - z^n = g(y^{n+1}) - (P_{\dt}g)(y^n) - \dt^{3/2} \sqrt{\frac{\gamma}{2\beta}} \nabla_p (\mL_{\gamma} g)(y^n)^{\T}  G^n\,.
\end{equation}
In view of the expression~\eqref{eqn:Lgamma} of the generator, we can express the last term in \eqref{eq:incrementZ}, up to a multiplicative factor~$\sqrt{\gamma/(2\beta)}$ as
\begin{equation}\label{eqn:part-2}
\begin{split}
& \nabla_q g(y^n)^{\T} M^{-1} G^n - \gamma \nabla_p g(y^n)^{\T} M^{-1} G^n + D_{qp}^2 g(y^n): (M^{-1}p^n) \otimes G^n\\
& \quad - D_{pp}^2 g(y^n) : (\nabla V(q^n) + \gamma M^{-1} p^n) \otimes G^n + \frac{\gamma}{\beta} \nabla_p \Delta_p g(y^n)^{\T} G^n\,.
\end{split}
\end{equation}
Subtracting the latter term (multiplied by~$\sqrt{\gamma/(2\beta)}$) from the order $\mathcal{O}(\dt^{3/2})$ term of~\eqref{eqn:part-1} leads to
\begin{equation}\label{eqn:order-3-half}
\begin{split}
\sqrt{\frac{\gamma}{2\beta}} D_{qp}^2 g(y^n) : (M^{-1}p^n) \otimes G^n - \sqrt{\frac{\gamma}{2\beta}} D_{pp}^2 g(y^n) : (\nabla V(q^n) + \gamma M^{-1} p^n) \otimes G^n\\
+ \frac{1}{6} \left(\frac{2\gamma}{\beta}\right)^{3/2} D_{ppp}^3 g(y^n) : (G^n)^{\otimes 3} - \frac{1}{\sqrt{2}} \left(\frac{\gamma}{\beta}\right)^{3/2} \nabla_p \Delta_p g(y^n)^{\T} G^n \PERIOD
\end{split}
\end{equation}
This provides the following expression for the martingale increment:
\begin{equation*}
\begin{split}
z^{n+1} - z^n &= \dt^{1/2} \sqrt{\frac{2\gamma}{\beta}}  \nabla_p g(y^n)^{\T} G^n 
+ \frac{\gamma \dt}{\beta} \left[ D_{pp}^2 g(y^n) : (G^n)^{\otimes 2} -\Delta_p g(y^n)\right]\\
&+ \dt^{3/2} \sqrt{\frac{\gamma}{2\beta}} \left[ D_{qp}^2 g(y^n) : (M^{-1}p^n) \otimes G^n
- D_{pp}^2 g(y^n) : (\nabla V(q^n) + \gamma M^{-1} p^n) \otimes G^n \phantom{\frac\gamma\beta}\right.\\
& \qquad \qquad \qquad + \left. \frac{2\gamma}{3\beta} D_{ppp}^3 g(y^n) : (G^n)^{\otimes 3}
- \frac{\gamma}{\beta} (\nabla_p \Delta_p g(y^n))^{\T} G^n\right]\,\PERIOD
\end{split}
\end{equation*}
It remains at this stage to find the function $g$. For a general perturbation $F(q,p)$ 
a possible solution is provided by Remark~\ref{rem:gfun}:
\[
g(q,p) = \frac{\beta}{2\gamma} \left(\nabla_p^* \nabla_p \right)^{-1} \nabla_p^* F(q,p)\,.
\]
A closed form expression for $g$ is thus available only in special cases, for example if the perturbing field is such that $F(q,p) = \nabla_p \phi(q,p)$ for some function $\phi$. In this case, we can consider $g(y) = \beta \phi(y)/(2\gamma)$ which leads to 
\begin{equation}\label{eqn:zinc}
\begin{split}
z^{n+1} &- z^n= \dt^{1/2}\sqrt{\frac{\beta}{2\gamma}} F(y^n)^{\T} G^n \\ 
              & \quad + \frac{\dt}{2} \left(\nabla_p F(y^n) : (G^n)^{\otimes 2} - (\mathrm{div}_p F)(y^n)\right)\\
              & \quad + \frac{\dt^{3/2}}{2}\sqrt{\frac{\beta}{2\gamma}}\left( (M^{-1}p^n)^{\T} \nabla_{q} F(y^n) G^n - (\nabla V(q^n) + \gamma M^{-1} p^n)^{\T} \nabla_p F(y^n) G^n \right)\\
              & \quad + \frac{\dt^{3/2}}{2}\sqrt{\frac{2\gamma}{\beta}}\left( \frac{1}{3} (D_{pp}^2 F)(y^n) : (G^n)^{\otimes 3}-\frac{1}{2} (\nabla_p \mathrm{div}_p F)(y^n)^{\T} G^n\right)\,.
\end{split}
\end{equation}
Note that the $\BIGO(\dt)$ term has mean zero since $\mathbb{E}[G_iG_j] = \delta_{ij}$.
Note that when $F$ depends on $q$ only, the above formula reduces to the formula~\eqref{eqn:2nd-order-weight} since, in this case, $\nabla_q g(q,p) = \beta/(2\gamma) \nabla F(q)^{\T} p$.

Let us finally provide the expression of the martingale increment 
$z^{n+1} - z^n$ for the scheme CBABC, for the specific case when 
$F$ only depends on $q$ and hence the function $g$ is 
of the form \eqref{eq:solutio_g_MP_Langevin}.
Recall the expansion~\eqref{eq:expansion_CBABC} associated with the scheme CBABC. 
Note first that the scheme CBABC involves two Gaussians $G_1^n$ and $G_2^n$
and the term of order $\dt^{1/2}$ in $z^{n+1} - z^n$ becomes
\[
\sqrt{\frac{\gamma}{\beta}} \nabla_p g(y^n)^{\T} (G_1^n + G_2^n),
\]
which can be seen as replacing the single Gaussian $G^n$ in~\eqref{eqn:part-1}
by $(G_1^n + G_2^n)/\sqrt{2}$.
Furthermore, the term of order $\dt$ vanishes when $F = F(q)$.
Finally, to obtain the term of order $\dt^{3/2}$, 
we start by noting that the last term in~\eqref{eqn:part-1} should be replaced with
\[
\sqrt{\frac{\gamma}{\beta}}\nabla_q g(y^n)^{\T} M^{-1}G_1^n -\frac14 \sqrt{\frac{\gamma}{\beta}} \gamma \nabla_p g(y^n)^{\T}  M^{-1} (3G_1^n+G_2^n) +\sqrt{\frac{\gamma}{\beta}} D_{qp}^2 g(y^n) : (M^{-1} p^n) \otimes (G_1^n+G_2^n).
\]
One should then subtract~\eqref{eqn:part-2} multiplied by a factor~$\sqrt{\gamma \beta^{-1}}/2$ and with~$G^n$ replaced by~$G_1^n+G_2^n$. Equation~\eqref{eqn:order-3-half} is then replaced by
\[
\frac{1}{2}\sqrt{\frac{\gamma}{\beta}} D_{qp}^2 g(y^n) : (M^{-1}p^n) \otimes (G_1^n+G_2^n) +
\frac12 \sqrt{\frac{\gamma}{\beta}} \left( \nabla_q g(y^n)-\frac{\gamma}{2}\nabla_p g(y^n) \right)^\T M^{-1}(G_1^n-G_2^n)
\, \PERIOD
\]
Plugging in $g$ of the form~\eqref{eq:solutio_g_MP_Langevin} leads to the formula~\eqref{eqn:2nd-order-weight-CBABC}.

\subsection*{Acknowledgements}
The work of P.P. was supported in part by the U.S. Army Research Office Award W911NF-19-1-0243.
The work of G.S. was partially funded by the European Research Council (ERC) under the European Union's Horizon 2020 research and innovation programme (grant agreement No 810367), and by the Agence Nationale de la Recherche under grants ANR-19-CE40-0010-01 (QuAMProcs) and ANR-21-CE40-0006 (SINEQ).
The research of T.W. was sponsored by the DEVCOM Army Research Laboratory and was accomplished under Cooperative Agreement Number W911NF-16-2-0190. 
The views and conclusions contained in this document are those of the authors and should not be interpreted as representing the official policies, either expressed or implied, of the Army Research Laboratory or the U.S. Government. The U.S. Government is authorized to reproduce and distribute reprints for Government purposes notwithstanding any copyright notation herein.

\end{document}